\documentclass[12pt,reqno,showkeys]{amsart}  
\usepackage{amssymb} 
\usepackage{amsthm}
\usepackage{amsmath}
\usepackage{mathrsfs}
\usepackage{color}
\usepackage{bbm}
\usepackage{verbatim}
\usepackage[utf8]{inputenc}
\usepackage[T1]{fontenc}
\usepackage[colorlinks=true, pdfstartview=FitV, linkcolor=blue,citecolor=blue, urlcolor=blue,dvipdfmx]{hyperref}
\allowdisplaybreaks[1]
\newtheorem{thm}{Theorem}[section]
\newtheorem*{def*}{Definition}

\newtheorem{prop}[thm]{Proposition}
\newtheorem{lem}[thm]{Lemma}
\newtheorem{rem}{Remark}[section]

\numberwithin{equation}{section}

\renewcommand{\theequation}{\thesection.\arabic{equation}}

\newenvironment{pr}[1]
   {{\noindent \bf Proof of {#1}.\  }}{\hfill \qed}
                
\newcommand{\pt}{\partial}     
             
\renewcommand{\th}{\theta}                
\newcommand{\D}{\mathcal{D}}
           
\newcommand{\F}{\mathcal {F}}

\newcommand{\R}{\mathbb R}

\newcommand{\Nt}{\mathbb N}
\newcommand{\N}{\nabla}

\newcommand{\al}{\alpha}

\newcommand{\gm}{\gamma}

\renewcommand{\L}{\mathcal{L}}

\newcommand{\ep}{\varepsilon}
\newcommand{\lm}{\lambda}

\newcommand{\Del}{\Delta}

\renewcommand{\t}{\tau}

\newcommand{\dx}{\,\mathrm{d}x}

\newcommand{\dsp}{\displaystyle}

\newcommand{\bt}{\beta}
\newcommand{\cd}{\cdot}

\renewcommand{\d}{\mathrm{d}}

\def\<{\langle }

\newcommand{\supp}{\text{supp~}}

\renewcommand{\qed}{\qquad\kern1pt   
   \vbox{\hrule height 0.6pt      
         \hbox{\vrule width 0.6pt 
               \vbox{\vskip 6pt}  
               \hskip 3pt
              \vrule width 1.3pt} 
         \hrule depth 1.3pt}     
   \kern1pt}

\newcommand{\eqntag}{\addtocounter{equation}{1}\tag{\theequation}} 

\begin{document}
\title[Global existence of solutions with critical mass]
{Global existence for the fully parabolic Keller--Segel system with critical mass on the plane}
\author{Tatsuya Hosono}
\address{Osaka Central Advanced Mathematical Institute, Osaka Metropolitan University,
	Osaka 558-8585, Japan 
	\& Laboratoire de Math\'{e}matiques (LAMA) UMR 5127,
	Universit\'{e} Savoie Mont Blanc, F--73000 Chamb\'{e}ry, France}
\email{tatsuya.hosono@omu.ac.jp}
\date{}
\pagestyle{myheadings}
\maketitle
\begin{abstract}
We study the global existence of solutions to the Cauchy problem for the
two-dimensional fully parabolic Keller--Segel system at the critical mass.
It is known that global-in-time existence holds for initial data with critical
mass under radial symmetry or suitable moment conditions, whereas the behavior
of general solutions in the critical regime remains delicate.
In this paper, we establish global-in-time existence for general initial data
with critical mass, without imposing any symmetry or moment assumptions.
The proof relies on the construction of a reconstructed Lyapunov functional,
combined with refined regularity estimates for the associated dissipative terms,
which enable us to control the solution dynamics in the critical regime.

\end{abstract}

\vspace{5mm}

\vspace{5mm}


\section{Introduction}
We study  the Cauchy problem for the parabolic-parabolic Keller--Segel system in~$\R^2$
\begin{equation}
\left\{
\begin{aligned}
&\pt_t u =\Delta u- \N\cd \left( u \N v \right),
& t>0,\,~ &x\in\R^2,
\\
&\pt_t v = \Delta v - \lambda v +u,
&t>0,\, ~&x\in\R^2,
\\
&(u,v)(0,x)=(u_{0},v_0)(x),
&\, &x\in\R^2
\end{aligned}
\right.
\label{eqn;KS}
\end{equation}
with a constant $\lambda\ge0$, where $u_0,v_0\ge 0$ on $\R^2$ and
$u_0,v_0\not\equiv0$.
The Keller--Segel system~\eqref{eqn;KS} is a fundamental mathematical model
of chemotaxis, describing chemotactic aggregation in the
cellular slime mold {\it Dictyostelium discoideum} during its life cycle \cite{KeSe,Pa}.
In this model, $u=u(t,x)$ and $v=v(t,x)$ denote the densities of cells and
the chemoattractant, respectively.
The parameter $\lambda$ represents the degradation rate of the chemical.
Cells migrate toward regions of higher concentrations of a chemical substance
secreted by the cells themselves.
From a biological viewpoint, it is natural to assume that the initial data
are nonnegative functions.
From a mathematical viewpoint, cell aggregation is interpreted as the
blowup of solutions at $t=T$ in the sense that
$\lim_{t\to T}\|u(t)\|_\infty=\infty.$

One of the central mathematical features of such systems is the phenomenon of critical
mass.
For positive sufficiently regular solutions $(u,v)$ to~\eqref{eqn;KS}, the first
component $u$ satisfies the mass conservation law 
$$\|u(t)\|_1=\|u_0\|_1\quad \text{for}~~~t>0.$$
Moreover, in two space dimensions, the global behavior of solutions is governed
by the size of the initial mass of $u$.
In particular, there exists a threshold value such that solutions with initial
mass below this threshold exist globally in time, whereas solutions with
initial mass above the threshold may blow up in finite time.
This critical-mass phenomenon has attracted considerable attention in the
literature.

The aim of this paper is to establish the global behavior of solutions to
\eqref{eqn;KS} whose initial mass is exactly equal to this threshold value.

\vspace{3mm}
The second equation in \eqref{eqn;KS} takes into account that cells are producing the chemoattractant themselves 
while this is diffusing into the environment. 
Since the chemoattractant attains its equilibrium on a time scale much quicker
than that of the cells, the simplified parabolic-elliptic system has also been
investigated \cite{JaLu,Na95}:
\begin{align*}
\left\{
\begin{aligned}
\partial_t u &= \Delta u - \nabla\cdot (u\nabla v),
& t>0,\ & x\in\R^2,\\
0 &= \Delta v - \lambda v + u,
& t>0,\ & x\in\R^2,
\\
&u(0,x)=u_{0}(x),
&\, &x\in\R^2.
\end{aligned}
\right.
\eqntag
\label{eqn;PKS}
\end{align*}
The system \eqref{eqn;PKS} is also related to models of gravitational interaction
of particles \cite{Bi-Na94,Wo}.
In~\eqref{eqn;PKS}, the second equation can be written as
$v = (-\Delta+\lambda)^{-1}u.$
As a consequence, the system~\eqref{eqn;PKS} can be reduced to a single equation
for~$u$.
Owing to this elliptic structure, the parabolic-elliptic system~\eqref{eqn;PKS} is more amenable to analysis than the fully parabolic system~\eqref{eqn;KS}, and has therefore been extensively studied and is now well understood in many aspects.
In particular, numerous works have been devoted to the critical-mass
phenomenon.
Indeed, the solution with $\|u_0\|_1\le 8\pi$ exists globally in time~\cite{BiKaLaNa06,Bl-Ca-Ca,Bl-Ca-Ma,BlDoPe,DoPe,Lo-Na-Ya1,Lo-Na-Ya2,Na11,NaOg11,NaOg16,Na-Se,Wei},
while the solution may blow up in finite time if $\|u_0\|_1>8\pi$~\cite{DoPe,Ko-Su,Ku-Og03,Wei}.
For the corresponding Cauchy--Neumann problem in bounded domains,
supplemented with homogeneous Neumann boundary conditions, see \cite{BiKaLaNa06disc,Bi-Na94,Ga-Za98,Na95,Na01} for instance.

\vspace{3mm}
Unlike the parabolic-elliptic system \eqref{eqn;PKS}, the fully parabolic system
\eqref{eqn;KS} is a strongly coupled parabolic system, and many approaches
developed for \eqref{eqn;PKS} are no longer applicable.
The global existence of solutions with sub-critical mass
$\|u_0\|_{1}<8\pi$ has nevertheless been established by combining
Lyapunov functionals with the Trudinger--Moser type inequality and its optimal
constant.
More precisely, the following results are known:
\begin{enumerate}
	\item[(i)]
	If $\|u_0\|_{1}<8\pi$, then the corresponding solution to
	\eqref{eqn;KS} exists globally in time \cite{CaCo08,Mi13,NaOg11}.
	
	\item[(ii)]
	If $\|u_0\|_{1}>8\pi$ and $(u_0,v_0)$ is radially symmetric, then there
	exists a radially symmetric solution that blows up in finite time
	\cite{Mi20,Mi20-SIAM}.
	
	\item[(iii)]
	If $\|u_0\|_{1}=8\pi$ and $(u_0,v_0)$ is either radially symmetric or
	satisfies the additional moment condition
	$u_0\ln(1+|x|^2)\in L^1(\R^2)$, then the corresponding solution exists globally
	in time.
	In contrast, for general initial data, the global behavior in the critical case
	remains delicate, and the solution either exists globally or blows up on the plane \cite{CaCo08,Mi13}.
\end{enumerate}
As for positive forward self-similar solutions to \eqref{eqn;KS}, refer to \cite{BiCoDo}.
For the corresponding Cauchy--Neumann problem in bounded domains, we refer to
\cite{Bi98,He-Ve,HoWa,Na-Se-Yo}, for instance.

As mentioned above, although several partial results are available, no complete global existence result has been obtained
for solutions with critical mass and general initial data without any
symmetry or moment assumptions.
One of the main difficulties lies in controlling the behavior of solutions at~$|x|\to\infty$.
Moreover, in the critical mass case, no global existence result is currently
known  for the corresponding Cauchy--Neumann problem to~\eqref{eqn;KS} in bounded domains.

In this paper, we establish the global existence of solutions with critical mass
for general initial data $(u_0,v_0)$.
To this end, we first recall the definition of solutions.

\vspace{2mm}
\begin{def*}
Let $4/3<p<2$ and $q=p/(p-1)$.
Given $T\in (0,\infty]$ and  $(u_0,v_0) \in L^1(\R^2)\times \dot{H}^{1}(\R^2)$,
we say that a pair of functions $(u,v)$ defined on $[0,T)\times \R^2$ is a mild solution to~\eqref{eqn;KS} on $[0,T)$ if
\begin{enumerate}
	\item $u\in C\big([0,T); L^1(\R^2)\big)\cap C\left( (0,T); L^p(\mathbb{R}^2) \right)$;
	\item $v\in C\big([0,T); \dot{H}^{1}(\R^2)\big)\cap C\left( (0,T); \dot{W}^{1,q}(\mathbb{R}^2) \right)$;
	\item $\dsp\sup_{0<t<T} t^{1-\frac1p}\|u(t)\|_{p}+\dsp\sup_{0<t<T}t^{\frac12-\frac1q}\|\nabla v(t)\|_{q}$;
	\item $(u,v)$ satisfies the integral formulations
	\begin{align*}
	&u(t)=e^{t\Del}u_0-\int_0^t \N\cd e^{(t-s)\Del} \big(u(s)\N v(s)\big) \,\d s,
	\\
	&v(t)= e^{t(\Del-\lm)} v_0+\int_0^t e^{(t-s)(\Del-\lm)} u(s) \,\d s
	\end{align*}
	for $0<t<T$, where $e^{t\Del}$ is the heat semigroup given by
	\begin{align*}
	e^{t\Del} f(x) := G_t*f(x)=\int_{ \R^2}G_t(x-y) f(y)\,\d y,
	\quad
	G_t(x)=\frac{1}{4\pi t} \exp\left(-\frac{|x|^2}{4 t}\right)
	\end{align*}
	for $x\in\R^2$ and $f \in L^p(\R^2)$.
\end{enumerate}
\end{def*}
The Lebesgue spaces $L^p(\R^2)$ are equipped with the norm $\|\cdot\|_{L^p(\R^2)}$,
which is abbreviated as $\|\cdot\|_p$ when there is no ambiguity.
We define the homogeneous Sobolev spaces by
\[
\dot{W}^{1,p}(\R^2)
:= \overline{C_0^\infty(\R^2)}^{\|\nabla (\cdot)\|_{L^p(\R^2)}}
\]
and $\dot{H}^1(\R^2):=\dot{W}^{1,2}(\R^2)$ when $p=2$.
$L_+^1(\R^2):=\{f \in L^1(\R^2); f\ge0~~ \text{and}~~ f\not\equiv0~~ \text{on $\R^2$} \}$.

\vspace{5mm}
The starting point for our analysis is the local well-posedness for \eqref{eqn;KS} and the solution~$(u,v)$ based on the above definition satisfies~\eqref{eqn;KS} in a classical sense on $(0,T)$, see Proposition~\ref{prop;LWP} in Section~\ref{sect;LWP} below.

We are now in a position to state our main result in this paper.
Main result in this paper reads as follows.
\vspace{2mm}
\begin{thm}\label{thm;global}
For $(u_0,v_0)\in L_+^1(\R^2) \times L^1_+(\R^2)\cap \dot{H}^1(\R^2)$,
let $(u,v)$ be the solution to \eqref{eqn;KS} on $(0,T)\times\R^2$.
Suppose that $\|u_0\|_1=8\pi$.
Then, the solution to \eqref{eqn;KS} exists globally in time.
\end{thm}
\vspace{2mm}

\begin{rem}
In Theorem~\ref{thm;global}, the initial data are only required to satisfy
\[
(u_0,v_0)\in L_+^1(\R^2)\times\bigl(L_+^1(\R^2)\cap \dot H^1(\R^2)\bigr)\quad \text{with}~~\|u_0\|_1=8\pi.
\]
In particular, no additional symmetry or moment assumptions are imposed.
\end{rem}
\vspace{2mm}
As mentioned above, the global existence of solutions with critical mass
without any symmetry or moment assumptions remains delicate.
In~\cite[Theorem~1.2]{Mi13}, it is shown that, for general initial data with
critical mass, solutions either exist globally in time or blow up on the plane,
by means of a contradiction argument.
Indeed, although in the sub-critical case regularity estimates can be obtained
by combining (modified) Lyapunov functionals with the Trudinger--Moser
inequality, such estimates are no longer directly available in the critical case due to the lack of regularity of solutions,
see~\eqref{eqn;afford} below.
This difficulty prevents the direct extension of classical entropy methods.

To overcome this difficulty, we introduce a reconstructed Lyapunov functional.
The main novelty of this work lies in the construction of a refined Lyapunov
functional specifically adapted to the whole space setting, which allows us
to control the behavior of solutions at $|x|\to\infty$ and to derive
regularity estimates even in the critical mass regime, without
imposing any symmetry or moment assumptions on the initial data.
As a consequence, we establish global-in-time existence for general initial
data at the critical mass, a result that was previously out of reach by existing
methods, see Subsection~\ref{sect;modified-Lyapunov} for details.
Moreover, the present approach is expected to be applicable to a broad class
of chemotaxis  systems in the whole space setting.

\vspace{3mm}
The remainder of this paper is organized as follows.
In Section~\ref{sect;preliminary}, we recall several preliminary lemmas needed
to prove the main result.
Section~\ref{sect;LWP} is devoted to the local well-posedness of~\eqref{eqn;KS}, based on the above definition of solutions.
Finally, in Section~\ref{sect;global}, we present the proof of
Theorem~\ref{thm;global}.

\vspace{5mm}
\section{Preliminary}\label{sect;preliminary}
In this section, we collect some tools to show our main result.
\begin{lem}\label{lem;L2andL3}
	For any $f\in H^1(\R^2)$ and nonnegative cut-off function $\phi\in C^\infty_0(\R^2)$,
	\begin{align*}
	\int_{\R^2}|f|^2\phi\dx\le\,&2\left(\int_{ \{|f|>1\}\cap \{\supp\phi\} }|f|\dx\right)\left(\int_{ \{|f|>1\} }\frac{|\nabla f|^2}{1+|f|} \phi\dx\right)
	\\
	&+4\left(\int_{\R^2}|f \nabla \phi^{\frac12}|\dx\right)^2+4\left(\int_{\R^2}|f|\phi\dx\right)
	\end{align*}
	and for any $\ep>0$
	\begin{align*}
	\int_{\R^2}|f|^3\phi\dx
	\le\,&\ep \left(\int_{  \{\supp\phi\} }(1+|f|)\ln(1+|f|)\dx\right)\left(\int_{\R^2 }|\nabla f|^2\phi\dx\right)
	\\
	&+C\left(\int_{\R^2}|f^{\frac32} \nabla \phi^{\frac12}|\dx\right)^2+C_\ep\left(\int_{\R^2}|f|\phi\dx\right),
	\end{align*}
	where the constant $C_\ep\to \infty$ as $\ep\to0$.
\end{lem}
The proof of Lemma \ref{lem;L2andL3} can be found in \cite[Lemma 2.2]{NaOg16}.

\begin{lem}\label{lem;Lp_Lq_heat}
	For $1 \le q \le p \le \infty$,
	let $f\in L^q (\R^n)$ and let $\al$ be a multi-index.
	Then it follows that
	\begin{equation*}
	\| \pt_x^{\al} e^{t\Del} f \|_p  
	\le C t^{-\frac{n}{2} \left(\frac{1}{q} - \frac{1}{p} \right) - \frac{|\al|}{2}}   \|f\|_q
	\end{equation*}
	for all $t>0$.
\end{lem}

The proof is immediately obtained by use of Young's inequality and the 
convolution expression of the heat evolution
by the heat kernel, see for instance \cite{GiGiSa}.
\begin{lem}\label{lem;v-Lpbound}
For $v_0\in L^1(\R^2)\cap \dot H^1(\R^2)$ and $f\in L^\infty \left(0,\infty;  L^1(\R^2) \right)$,
let $v$ be the solution to $\partial_t v =\Delta v -\lambda v +f$  on $(0,\infty)\times\R^2$ with the initial data $v_0$ and $\lambda\ge0$. Then, for any $1\le p<\infty$, 
\begin{align*}
\| v(t)\|_p
\le\,&
\left\{
\begin{aligned}
&\|v_0\|_p+\frac1p\|f\|_{L^\infty\left(0,\infty; L^1(\R^2)\right)}  t^{\frac1p} &\text{if $\lambda=0$},
\\
&\|v_0\|_p+\lambda^{-\frac1p} \|f\|_{L^\infty\left(0,\infty; L^1(\R^2)\right)}  \Gamma\left(\frac1p\right)&\text{if $\lambda>0$}.
\end{aligned}
\right.
\end{align*}
When $p=1$,
\begin{align*}
\int_{\R^2 }v\,\d x
=
\left\{
\begin{aligned}
&\|v_0\|_1+ t \|u_0\|_1 &\text{if}~\lambda=0,
\\
&e^{-\lambda t}\|v_0\|_1+\frac{1-e^{-\lambda t}}{\lambda}\|u_0\|_1&\text{if}~\lambda>0.
\end{aligned}
\right.
\end{align*}
\end{lem}
The proof of Lemma \ref{lem;v-Lpbound} is obtained by the direct computations using the integral formulation for $v$.

\vspace{5mm}
\section{Local-in-time solutions}\label{sect;LWP}
This section deals with the local-in-time solution $(u,v)$ to \eqref{eqn;KS} corresponding to the initial data $(u_0,v_0)\in L^1(\R^2)\times \dot{H}^1(\R^2)$.
\begin{prop}[local-in-time solution]\label{prop;LWP}
	Let $4/3<p<2$ and $q=p/(p-1)$.
	For  $(u_0,v_0)\in L^1(\R^2)\times  \dot{H}^{1}(\R^2)$,
	there exist $T\in (0,\infty]$ and a unique mild solution $(u,v)$ to~\eqref{eqn;KS} on $[0,T)$ in the sense of the definition.
	Besides, the solution $(u,v)$ has higher regularity on $(0,T)$ and satisfies the system~\eqref{eqn;KS} in a classical sense on $(0,T)\times\R^2$.
	In addition,
	\begin{equation*}
	\text{either $T=\infty$ or $T<\infty$ and } \lim_{t\to T} \left\{ \|u(t)\|_{r} + \|\nabla v(t)\|_{r+1}\right\} = \infty,~~ 1<r\le \infty. 
	\eqntag
	\label{eqn;bucrit}
	\end{equation*}	
	Moreover, if $u_0,v_0 \ge 0$ then $u,v> 0$ in $(0,T)\times\R^2$ and
	\begin{equation*}
	\|u(t)\|_1=\|u_0\|_1, \qquad t\in [0,T).	
	\end{equation*} 
\end{prop}
Let us show the proof of Proposition \ref{prop;LWP}.
We fix constants $\eta>0$ and $T>0$ to be chosen later depending on the initial data.
Let $M=4(\|u_0\|_1+\|\nabla v_0\|_2)$ and let $4/3<p<2$ and $q=p/(p-1)$.
Define the space $\mathcal{X}_T$ by
\begin{align*}
\mathcal{X}_T:=
\left\{
\begin{aligned}
&u \in L^\infty \left( 0,T;L^1(\mathbb{R}^2) \right)\cap C\left( (0,T); L^p(\mathbb{R}^2)  \right),
\\
& v \in L^\infty\left( 0,T; \dot{H}^1(\mathbb{R}^2) \right)\cap C\left( (0,T); \dot{W}^{1,q}(\mathbb{R}^2) \right);
\\
&\sup_{0<t<T}\|u(t)\|_1+\sup_{0<t<T}\|\nabla v(t)\|_2 \le M,
\\
&\sup_{0<t<T}t^{1-\frac1p}\|u(t)\|_p\le \eta_1,~~\sup_{0<t<T}t^{\frac12-\frac1q}\|\nabla v(t)\|_q \le \eta_2
\end{aligned}
\right\}
\end{align*}
as well as the distance by
\begin{align*}
d_{\mathcal{X}_T} \left( (u,v), (\widetilde{u},\widetilde{v}) \right):=\sup_{0<t<T}t^{1-\frac1p}\|u(t)-\widetilde{u}(t)\|_p
+\sup_{0<t<T}t^{\frac12-\frac1q}\|\nabla v(t)-\widetilde{v}(t)\|_q.
\end{align*}
We then realize that $(\mathcal{X}_T, d_{\mathcal{X}_T})$ is a complete metric space.
Next, the mapping $\Phi=\left(\Phi_1, \Phi_2\right): \mathcal{X}_T \to \mathcal{X}_T$ is denoted by
\begin{align*}
&\Phi_1[u,v](t):=e^{t\Delta} u_0 - \int_0^t \nabla\cdot e^{(t-s)\Delta} \left(u(s)\nabla v(s)\right)\,\d s,
\\
&\Phi_2[u,v](t):=e^{t(\Delta-\lambda)} v_0+\int_0^t e^{(t-s)(\Delta-\lambda)} u(s)\,\d s
\end{align*}
for $(u,v)\in \mathcal{X}_T$ and $t\in(0,T)$.

We shall show that $\Phi$ is a contraction mapping from $X_T$ to itself by taking $(\eta_1,\eta_2,T)$ properly in order to apply the Banach fixed point theorem.

\begin{lem}\label{lem;fixed-point}
Let $(u,v)\in \mathcal{X}_T$.
Then, there exist constants $(\eta_1,\eta_2,T)>0$ sufficiently small such that
\begin{align*}
 \sup_{0<t<T}\|\Phi_1[u,v](t)\|_1+\sup_{0<t<T}\|\nabla \Phi_2[u,v](t)\|_2\le\, M
\end{align*}
and
\begin{align*}
\sup_{0<t<T}t^{1-\frac1p}\|\Phi_1[u,v](t)\|_p\le\,\eta_1,~~
\sup_{0<t<T}t^{\frac12-\frac1q}\|\nabla \Phi_2[u,v](t)\|_q\le\,\eta_2.
\end{align*}
Besides, for $(u,v)\in \mathcal{X}_T$ and $(\widetilde{u},\widetilde{v})\in \mathcal{X}_T$,
\begin{align*}
d_{\mathcal{X}_T}\left(\Phi[u,v],\Phi[\widetilde{u},\widetilde{v}]\right)
\le \frac12 d_{\mathcal{X}_T}\left((u,v),(\widetilde{u},\widetilde{v})\right).
\eqntag\label{eqn;contraction-d}
\end{align*}
\end{lem}
\begin{pr}{Lemma \ref{lem;fixed-point}}
Let $(u,v)\in \mathcal{X}_T$. Owing to H\"{o}lder's inequality and the properties of the heat kernel,
\begin{align*}
&\|\Phi_1[u,v](t)\|_1
\\
\le\,& \| e^{t\Del} u_0\|_1
+C\int_0^t (t-s)^{-\frac12} \|u(s)\N v(s)\|_1\, \d s
\\
\le\,& \|  u_0\|_1
+C\int_0^t (t-s)^{-\frac12}\|u(s)\|_{p}\|\N v(s)\|_q \,\d s
\\
\le\,&\|u_0\|_1+C \int_0^t (t-s)^{-\frac12} s^{-\frac32+\frac1p+\frac1q}\,\d s
\left(\sup_{0< s<t}s^{1-\frac1p}\|u(s)\|_{p}\right)
\left(\sup_{0< s<t} s^{\frac12-\frac1q}\|\nabla v(s)\|_{q}\right)
\\
\le\,&\frac{M}{4}+C_1 \mathsf{B}\left(\frac12,\frac12\right)  \eta_1\eta_2,
\eqntag \label{eqn;L1}
\end{align*}
where 
\begin{align*}
\frac1p+\frac1q=1
\end{align*}
and the beta function $\mathsf{B}(x,y)$ for $x,y>0$ is given by
\begin{align*}
\mathsf{B}(x,y):=\int_0^1(1-\t)^{x-1} \t^{y-1}\,\d\t,
\quad
(x,y)\in(0,\infty)^2.
\end{align*}
Choosing $\eta_1>0$ small such that
\begin{equation*}
C_1\mathsf{B}\left(\frac12,\frac12\right)  \eta_1\eta_2\le \frac M4,
\end{equation*}
we see that
\begin{align*}
\sup_{0< t<T}\|\Phi_1[u,v](t)\|_1 \le \frac{M}{2}.
\end{align*}
For the estimate for $\Phi_2[u,v]$, the smoothing properties of the heat kernel imply that
\begin{align*}
\|\nabla \Phi_2[u,v](t)\|_{2}
\le\,&\|\nabla v_0\|_{2}+C\int_0^t (t-s)^{-\left(\frac1p-\frac12\right)-\frac12} s^{-1+\frac1p}\,\d s
\left(\sup_{0< s<t}s^{1-\frac1p}\|u(s)\|_{p}\right)
\\
\le\,& \frac{M}{4} + C_2 \mathsf{B}\left(1-\frac1p, \frac1p\right) \eta_1.
\eqntag\label{eqn;W-22}
\end{align*}
Taking $\eta_1>0$ such that
\begin{align*}
C_2 \mathsf{B}\left(1-\frac1p, \frac1p\right) \eta_1\le \frac{M}{4},
\end{align*}
it follows that
\begin{align*}
\sup_{0< t<T}\|\nabla \Phi_2[u,v](t)\|_{2} \le \frac{M}{2}.
\end{align*}
Next 
 by H\"{o}lder's inequality,
\begin{align*}
&t^{1-\frac1p}\left\|\int_0^t \N\cdot e^{(t-s)\Del}(u(s)\N v(s))ds\right\|_{p}
\\
\le\,&Ct^{1-\frac1p}
\int_0^t (t-s)^{-\left(1-\frac1p\right)-\frac12} \| u(s)\N v(s)\|_1 \,\d s
\\
\le\,&Ct^{1-\frac1p}
\int_0^t (t-s)^{\frac1{p}-\frac32} \|u(s)\|_{p} \|\N v(s)\|_{q}\, \d s
\\
\le\,&Ct^{1-\frac1p}
\int_0^t (t-s)^{\frac1{p}-\frac32} s^{-\frac32+\frac1p+\frac1q}\,\d s
\left(\sup_{0< s<t}s^{1-\frac1p}\|u(s)\|_{p}\right)\left(\sup_{0< s<t}s^{\frac12-\frac1q}\|\nabla v(s)\|_{q}\right)
\\
\le\,&C_4 \mathsf{B}\left(\frac1p-\frac12,\frac12\right) \eta_2 \left(\sup_{0< s<t}s^{1-\frac1p}\|u(s)\|_{p}\right).
\eqntag \label{eqn;L-4/3}
\end{align*}
Choosing
\begin{equation*}
C_4 \mathsf{B}\left(\frac1p-\frac12,\frac12\right) \eta_2\le \frac12
\end{equation*}
in~\eqref{eqn;L-4/3} gives
\begin{align*}
t^{1-\frac1p}\|\Phi_1[u,v](t)\|_{p}\le
t^{1-\frac1p}\|e^{t\Del}u_0\|_{p} + \frac12 \left(\sup_{0< s<t}s^{1-\frac1p}\|u(s)\|_{p}\right),
\end{align*}
from which we deduce that 
\begin{align*}
\sup_{0< t<\tau}t^{1-\frac1p}\|\Phi_1[u,v](t)\|_{p}\le
\sup_{0< t<\tau}t^{1-\frac1p}\|e^{t\Del}u_0\|_{p}+\frac12 \left(\sup_{0< t<\tau}t^{1-\frac1p}\|u(t)\|_{p}\right)
\eqntag \label{eqn;L4/3}
\end{align*}
for all $\tau\in (0,T]$. In particular, since
\begin{align*}
\lim_{t\to 0}t^{1-\frac1p}\|e^{t\Del}u_0\|_{p}=0
\eqntag
\label{eqn;Lpsmall}
\end{align*}
according to Weissler \cite[equation~(3.4)]{Weis},
we can choose $T$ small enough such that
\begin{align*}
t^{1-\frac1p}\|e^{t\Del}u_0\|_{p}\le \frac{\eta_1}{2},\quad
t\in(0,T).
\eqntag
\label{eqn;initial-0}
\end{align*}
Combining~\eqref{eqn;L4/3} and~\eqref{eqn;initial-0} leads us to
\begin{equation*}
\sup_{0< t<T}t^{1-\frac1p}\|\Phi_1[u,v](t)\|_{p}\le \frac{\eta_1}{2}+\frac{\eta_1}{2}=\eta_1.
\end{equation*}
As for the estimates of $\Phi_2[u,v]$, the proof is similar to that of~\eqref{eqn;W-22}.
Indeed, due to the linearity of $\Phi_2$
\begin{align*}
t^{\frac12-\frac1q}\|\nabla \Phi_2[u,v](t)\|_{q}\le 
t^{\frac12-\frac1q}\|\nabla e^{t\Del} v_0\|_{q}
+ C_5 \mathsf{B}\left(\frac12-\frac1p+\frac1q,\frac1p\right) \eta_1,
\eqntag \label{eqn;W-2p}
\end{align*}
observing that the right-hand side of~\eqref{eqn;W-2p} is finite due to the constraint $q<2p/(2-p)$. Since
\begin{equation*}
\lim_{t\to 0}t^{\frac12-\frac1q}\|\nabla e^{t\Del} v_0\|_{q}=0,
\end{equation*}
we may choose $\eta_1,T>0$ small enough so that 
\begin{align*}
C_5 \mathsf{B}\left(\frac12-\frac1p+\frac1q,\frac1p\right) \eta_1
\le \frac{\eta_2}{2}, 
\end{align*}
and 
\begin{equation*}
\sup_{0< t<T}t^{\frac12-\frac1q}\|\nabla e^{t\Del} v_0\|_{q} \le \frac{\eta_2}{2}, 
\end{equation*}
and thereby deduce from~\eqref{eqn;W-2p} that 
\begin{equation*}
\sup_{0<t<T}t^{\frac12-\frac1q} \|\nabla \Phi_2[u,v](t)\|_{q}\le \eta_2.
\end{equation*}
Finally, an analogous argument implies the assertion~\eqref{eqn;contraction-d}.
\end{pr}

\vspace{5mm}
Lemma \ref{lem;fixed-point} gives that $\Phi:\mathcal{X}_T\to \mathcal{X}_T$ is a contraction mapping, so that
one may find a unique fixed point $(u,v)\in\mathcal{X}_T$ such that $\Phi_1[u,v]=u$ and $\Phi_2[u,v]=v$.
Next, we shall show the continuous dependence on the initial data.

\begin{lem}\label{lem;continuous-depend}
There exists a unique mild solution $(u,v)$ to \eqref{eqn;KS} defined on a maximal time interval $[0,T)$ with $T\in (0,\infty]$ which  depends continuously on the initial data. More precisely, given $(u_0,v_0)$ and $(\widetilde{u}_0, \widetilde{v}_0)$ in $L^1(\R^2)\times \dot{H}^{1}(\R^2)$, we denote by $(u,v)$ and $(\widetilde{u},\widetilde{v})$ the corresponding solutions to~\eqref{eqn;KS}, respectively, with respective maximal existence time $T$ and $\widetilde{T}$. For any $0<t_0<\min\{T,\widetilde{T}\}$, there is $C(t_0)>0$ such that
\begin{align*}
&\sup_{0< t<t_0}\| u(t)-\widetilde{u}(t)\|_1 \le C(t_0) \left(\|u_0-\widetilde{u}_0\|_1+\|\nabla (v_0-\widetilde{v}_0)\|_{2}\right),
\\
&\sup_{0< t<t_0}\| \nabla (v(t)-\widetilde{v}(t))\|_{2} \le C(t_0) \left(\|u_0-\widetilde{u}_0\|_1+\|\nabla (v_0-\widetilde{v}_0)\|_{2}\right).
\end{align*}
\end{lem}

\begin{pr}{Lemma \ref{lem;continuous-depend}}
Let us show that  $(u,v)$ established in Lemma \ref{lem;fixed-point} is continuous with the initial data.
Recalling \eqref{eqn;L4/3} with $\Phi_1[u,v]=u$, we have
\begin{align*}
\sup_{0< s<t}  t^{1-\frac1p}\|u(t)\|_{p}\le 2 \sup_{0< s<t} t^{1-\frac1p}\|e^{t\Del}u_0\|_{p}
\end{align*}
for $0<t<T$.
According to \eqref{eqn;Lpsmall}, for any $\ep>0$, there exists a constant $t_\ep\in(0,T)$ such that
\begin{align*}
\sup_{0< t<t_\ep}t^{1-\frac1p}\|e^{t\Del}u_0\|_{p} \le\,\frac{\ep}{2}.
\end{align*}
Together with the above estimate, we obtain
\begin{align*}
\sup_{0< t<t_\ep}  t^{1-\frac1p}\|u(t)\|_{p}\le \ep,
\end{align*}
that is,
\begin{align*}
\lim_{t\to 0}\sup_{0< s<t}  t^{1-\frac1p}\|u(t)\|_{p}=0.
\eqntag
\label{eqn;Lpsmall_2}
\end{align*}
By the similar argument to that of \eqref{eqn;L1},
\begin{align*}
\|u(t)-u_0\|_1\le\,&\| e^{t\Delta} u_0 -u_0\|_1+C_1 \mathsf{B}\left(\frac12,\frac12\right)  \left(\sup_{0< s<t}s^{1-\frac1p}\|u(s)\|_{p}\right)
\left(\sup_{0< s<t} s^{\frac12-\frac1q}\|\nabla v(s)\|_{q}\right),
\end{align*}
from which it follows from the continuity of the heat kernel that
\begin{align*}
\lim_{t\to 0}\|u(t)-u_0\|_1=0,
\end{align*}
which means that $u\in C\left( [0,T); L^1(\mathbb{R}^2) \right)$.
Analogously, from \eqref{eqn;W-2p} with $\Phi_2[u,v]=v$ and \eqref{eqn;Lpsmall_2} 
\begin{align*}
\lim_{t\to 0}\sup_{0< s<t}s^{\frac12-\frac1q}\|\nabla v(s)\|_{q}=0,
\end{align*}
so that $v\in C\left( [0,T); \dot{H}^1(\mathbb{R}^2) \right)$, thereby proving the mild solution $(u,v)$ to \eqref{eqn;KS} on $(0,T)$.

Next,  show the uniqueness of solutions. Let $(u,v)$ and $(\widetilde{u},\widetilde{v})$ be two solutions to~\eqref{eqn;KS} on~$[0,T)$ corresponding to the initial data $(u_0,v_0)$. Noting that the mild formulation provides the representation formula for $v$
\begin{equation*}
v(t)-\widetilde{v}(t)=\int_0^t e^{(t-s)(\Del-\lm)} (u(s)-\widetilde{u}(s))ds,
\end{equation*}
we obtain by the similar way to that of \eqref{eqn;W-2p}
\begin{align*}
\sup_{0< s<t}s^{\frac12-\frac1q}\| \nabla (v(s) - \widetilde{v} (s)) \|_{q} 
\le\,& C \sup_{0< s<t} s^{1-\frac1p} \| u(s) -\widetilde{u}(s) \|_{p}.
\end{align*}
Using this estimate implies from the representation formula for $u-\widetilde{u}$ that
\begin{align*}
&\sup_{0< s<t}s^{1-\frac1p} \| u(s) - \widetilde{u}(s)\|_p
\\
\le\,&C \left(\sup_{0< s<t}s^{\frac12-\frac1q}\|\nabla v(s)\|_{q}\right)
\left(\sup_{0< s<t}s^{1-\frac1p}\|u(s)-\widetilde{u}(s)\|_{p}\right)
\\
&+ C\left(\sup_{0< s<t}s^{1-\frac1p}\|u(s)\|_{p}\right)\left(\sup_{0< s<t}s^{\frac12-\frac1q}\|\nabla (v(s)-\widetilde{v}(s))\|_{q}\right)
\\
\le\,&C \left( \sup_{0< s<t}s^{1-\frac1p}\|u(s)\|_{p}+\sup_{0< s<t}s^{\frac12-\frac1q}\|\nabla v(s)\|_{q}\right)
\left(\sup_{0< s<t}s^{1-\frac1p}\|u(s)-\widetilde{u}(s)\|_{p}\right).
\end{align*}
Since there exists $t_0\in(0,T)$ sufficiently small such that
\begin{align*}
C \left( \sup_{0< s<t_0}s^{1-\frac1p}\|u(s)\|_{p}+\sup_{0< s<t_0}s^{\frac12-\frac1q}\|\nabla v(s)\|_{q}\right)\le\, \frac12,
\end{align*}
we have
\begin{align*}
\sup_{0< t<t_0} t^{1-\frac1p} \| u(t)-\widetilde{u}(t)\|_{p}=0,
\end{align*}
and hence, we conclude that $u(t)=\widetilde{u}(t)$ for $0<t<t_0$, which in turn implies that $v(t)=\widetilde{v}(t)$ for $0<t<t_0$. Define $t_1\ge 0$ by
\begin{equation*}
t_1 := \sup\left\{ \tau>0\,;\,(u,v)(t)=(\widetilde{u},\widetilde{v})(t),\,\text{for all  } 0\le t<\tau \right\}.
\end{equation*}
Clearly $0<t_0\le t_1\le T$. Let us then suppose for contradiction that $t_1<T$. Repeating the above argument with initial data $(u,v)(t_1)=(\widetilde{u},\widetilde{v})(t_1)$, there exists $t_2\in(0,T-t_1)$ such that $(u,v)(t)=(\widetilde{u},\widetilde{v})(t)$ for all $0 \le t \le t_1+t_2$, which contradicts the definition of~$t_1$. 
Consequently, $(u,v)(t)=(\widetilde{u},\widetilde{v})(t)$ for $0<t<T$ and we have established the claimed uniqueness.

Given for $(u_0,v_0)$ and $(\widetilde{u}_0, \widetilde{v}_0)$ in $L^1(\R^2)\times \dot{H}^{1}(\R^2)$, we denote by $(u,v)$ and $(\widetilde{u},\widetilde{v})$ the corresponding solutions to~\eqref{eqn;KS}, respectively, with respective maximal existence time $T$ and~$\widetilde{T}$.
For $0<t_0<\min\{T,\bar{T}\}$ and $0<t<t_0$, a similar argument as above follows that
\begin{align*}
&\sup_{0< s<t} s^{1-\frac1p} \| u(s) - \widetilde{u}(s)\|_{p} 
\\
\le\,& C \| u_0 - \widetilde{u}_0\|_1+C\left(\sup_{0< s<t}s^{\frac12-\frac1q}\|\nabla v(s)\|_{q}\right)
\left(\sup_{0< s<t}s^{1-\frac1p}\|u(s)-\widetilde{u}(s)\|_{p}\right)
\\
&+C \left(\sup_{0< s<t}s^{1-\frac1p}\|u(s)\|_{p}\right)\left(\sup_{0< s<t}s^{\frac12-\frac1q}\|\nabla (v(s)-\widetilde{v}(s))\|_{q}\right)
\\
\le\,&C \| u_0 - \widetilde{u}_0\|_1+C\left(\sup_{0< s<t}s^{\frac12-\frac1q}\|\nabla v(s)\|_{q}\right)
\left(\sup_{0< s<t}s^{1-\frac1p}\|u(s)-\widetilde{u}(s)\|_{p}\right)
\\
&+C \left(\sup_{0< s<t}s^{1-\frac1p}\|u(s)\|_{p}\right)
\left(\| \nabla (v_0-\widetilde{v}_0)\|_{2}+\sup_{0< s<t}s^{1-\frac1p}\|u(s)-\widetilde{u}(s)\|_{p}\right),
\end{align*}
which implies that
\begin{align*}
\sup_{0< s<t}s^{1-\frac1p}\|u(s)-\widetilde{u}(s)\|_{p}\le
C(t_0) \left(\|u_0-\widetilde{u}_0\|_1+\|\nabla (v_0-\widetilde{v}_0)\|_{2}\right).
\end{align*}
Moreover, we observe
\begin{align*}
\sup_{0< s<t}\| \nabla (v(t) - \widetilde{v}(s))\|_{2}
\le\,& C \| \nabla(v_0 - \widetilde{v}_0)\|_{2} + C \sup_{0< s<t}s^{1-\frac1p}\|u(s)-\widetilde{u}(s)\|_{p}
\\
\le\,&C(t_0) \left(\|u_0-\widetilde{u}_0\|_1+\|\nabla (v_0-\widetilde{v}_0)\|_{2}\right).
\end{align*}
Hence,
\begin{align*}
\sup_{0< s<t}\| u(s)-\widetilde{u}(s)\|_1 
\le\,&C(t_0) \left(\|u_0-\widetilde{u}_0\|_1+\|\nabla (v_0-\widetilde{v}_0)\|_{2}\right).
\end{align*}
As a consequence, the local well-posedness for~\eqref{eqn;KS} follows.
\end{pr}

\vspace{5mm}
\begin{pr}{Proposition \ref{prop;LWP}}
By virtue of Lemmas~\ref{lem;fixed-point} and~\ref{lem;continuous-depend},
we have established the existence of a unique mild solution to~\eqref{eqn;KS} defined on a maximal time interval $[0,T)$ with~$T\le\infty$. As for the regularity of the solutions, we use the standard iteration argument with respect to the derivative. Define $|\N|^\al f(x):=\F^{-1}[ |\xi|^\al \F f(\xi)\,](x)$  for $x\in\R^2$ and $\al>0$, where $\F$ denotes the Fourier transform. Let $\tau\in (0,T)$ and $t\in (0,\tau)$. Recalling that, for $4/3<p<2$,
\begin{equation*}
U_{p,0}(\tau) := \sup_{0< s<\tau}s^{1-\frac1p}\|u(s)\|_{p}<+\infty, \quad V_{q,0}(\tau) := \sup_{0< s<\tau}s^{\frac12-\frac1q}\|\nabla v(s)\|_{q}<+\infty,
\end{equation*}
by Lemma~\ref{lem;fixed-point}, we may argue as in the proof of~\eqref{eqn;L-4/3} to obtain
\begin{align*}
&t^{1-\frac1p+\frac{\al}{2}} \| |\N|^{\al} u(t) \|_{p}
\\
\le\,& t^{1-\frac1p+\frac{\al}{2}} \| |\N|^\al e^{t\Del}u_0\|_{p} + C U_{p,0}(\tau) V_{q,0}(\tau) t^{1-\frac1p+\frac{\al}{2}}  \int_0^t  (t-s)^{-1+\frac1{p}-\frac{1+\al}2} s^{-1+\frac{1}p} s^{-\frac12+\frac1q}\,\d s
\\
\le\,&C\|u_0\|_1 + C  U_{p,0}(\tau) V_{q,0}(\tau) \mathsf{B}\left(\frac1p-\frac{1+\al}2,\frac12\right),
\end{align*}
where $4/3<p<2$ and
\begin{align*}
\al <\frac2p-1<\frac12,
\end{align*}
which implies that
\begin{align*}
U_{p.\al}(\tau):=\sup_{0< t<\tau} t^{1-\frac1p+\frac{\al}{2}} \| |\N|^{\al} u(t)\|_{p} <+\infty,
\quad
\al<\frac2p-1.
\end{align*}
This property yields additional regularity estimates for $v$. Indeed, for $\bt\ge \al$, 
\begin{align*}
\| |\N|^{1+\beta} v(t) \|_{q}
\le\,&C t^{-\frac12+\frac1q-\frac\bt2} \| \nabla v_0\|_{2} + C U_{p,\al}(\tau) \int_0^t (t-s)^{-\left(\frac1p-\frac1q\right)-\frac{1+\bt-\al}{2} } s^{-1+\frac1p-\frac\al2}\,\d s 
\\
=\,&C  t^{-\frac12+\frac1q-\frac\bt2} \| \nabla v_0\|_{2} 
+C U_{p,\al}(\tau)  t^{-\frac12+\frac1q-\frac\bt2} \mathsf{B}\left(\frac12-\frac1p+\frac1q-\frac{\bt-\al}{2}, \frac1p-\frac{\al}{2}\right),
\end{align*}
provided that
\begin{align*}
\al\le \bt<\al+3-\frac4p<2-\frac2p<1,
\end{align*}
from which it follows that 
\begin{align*}
V_{q,\bt}(\tau) := \sup_{0< t<\tau} t^{\frac12-\frac1q+\frac\bt2} \| |\N|^{1+\bt} v(t) \|_{q}<+\infty,
\quad  \bt<\al+3-\frac4p\in (\al,1).
\end{align*}
The previous analysis leads us to a higher regularity estimate for $u$. Indeed, let $\th\ge\al$. Since $\al\le\bt\le\gm$,
\begin{align*}
& \| |\N|^{\th}u(t)\|_{p} \\
\le\,&C t^{-\left(1-\frac1p\right)-\frac\th2}\|u_0\|_1
\\
&+C\int_0^t(t-s)^{-\left(1-\frac1p\right)-\frac{1+\th-\al}{2}} \left(\||\N|^{\al} u(s)\|_{p}\|\nabla v(s)\|_q+\|u(s)\|_{p}\| |\N|^{1+\al} v(s)\|_{q}\right)\,\d s
\\
\le\,&Ct^{-\left(1-\frac1p\right)-\frac\th2}\|u_0\|_1
\\
&+ C \big(U_{p,\alpha}(\tau) V_{q,0}(\tau) + U_{p,0}(\tau) V_{q,\al}(\tau)\big) 
\int_0^t(t-s)^{-\left(1-\frac1p\right)-\frac{1+\th-\al}{2}} s^{-\frac32+\frac1p+\frac1q-\frac{\al}{2}}\,\d s
\\
=\,&Ct^{-\left(1-\frac1p\right)-\frac\th2}\|u_0\|_1
\\
& + C \big(U_{p,\al}(\tau) V_{q,0}(\tau) + U_{p,\al}(\tau) V_{q,0}(\tau)\big) t^{-\left(1-\frac1p\right)-\frac\th2} \mathsf{B}\left(\frac1p-\frac{1+\th-\al}{2},\frac12-\frac\al2\right),
\end{align*}
where
\begin{align*}
\al \le \th<\al+\frac2p-1<1,\quad
\al<\frac2p-1.
\end{align*}
Hence, 
\begin{align*}
U_{p,\th}(\tau) = \sup_{0< t< \tau}t^{1-\frac1p+\frac{\th}{2}}\||\N|^\th u(t)\|_{p}<+\infty,
\quad
\al \le \th<\al+\frac2p-1<1,\quad
\al<\frac2p-1.
\end{align*}
Iterating further the above arguments, we obtain that the solution possesses higher regularity. In fact, using the uniqueness result established in Lemma~\ref{lem;continuous-depend}, we may solve the problem~\eqref{eqn;KS} with the initial data $(u(t_0),v(t_0))$ for some $t_0\in (0,T)$. We thus eventually show that for $s>0$, $u(t)\in L^1(\R^2)\cap W^{s,p}(\R^2)$,
$v(t)\in  \dot{H}^{1}(\R^2)\cap W^{s,q}(\R^2)$
for $t_0<t<T$, 
and satisfies the system~\eqref{eqn;KS} in a  classical sense on $(0,T)\times\R^2$.
Also, classical arguments based on the proof of Lemma~\ref{lem;fixed-point} lead to the blowup criterion~\eqref{eqn;bucrit}.

Let  $(u_0,v_0)$ be nonnegative initial conditions with $u_0,v_0\not\equiv0$ and
$\{u_{0,k}\}_{k\in\Nt }$ and
$\{v_{0,k}\}_{k\in\Nt }$ be sequences
of nonnegative functions in $C_0^{\infty}(\R^2)$ converging to $u_0$ in $L^{1}(\R^2)$ and $v_0$ in $\dot{H}^{1}(\R^2)$  as $k\to\infty$, respectively. Then, for each $k\in\Nt$, we can construct the unique, smooth and integrable solution
to~\eqref{eqn;KS} corresponding to the initial data $(u_{0,k},v_{0,k})$.
By the parabolic regularity theory,
the solutions $(u_k,v_k)$ belong to $C^{\infty}( [0,T)\times \R^2)$
and are positive on $(0,T)\times\R^2$. In addition, by the continuous dependence on the initial data established in Lemma~\ref{lem;continuous-depend}, we have the convergences $u_k\to u$ in $C\big([0,T); L^1(\R^2)\big)$ and $v_k\to v$ in $C\big([0,T); \dot{H}^{1}(\R^2)\big)$ as $ k\to \infty$,
where $(u,v)$ is the solution to~\eqref{eqn;KS} with initial data $(u_0,v_0)$. Hence, one can find subsequences of the approximated solutions
which converge to $(u,v)$ almost everywhere in $(0,T)\times\R^2$, which implies that $u$ and $v$ are positive on $(0,T)\times\R^2$.
Moreover the mass conservation law also holds by integrating $u$ over $\R^2$.
\end{pr}

\vspace{5mm}
\section{Global existence of solutions with critical mass}\label{sect;global}

In this section, we prove the global existence of solutions to \eqref{eqn;KS}
with critical mass, as stated in Theorem~\ref{thm;global}.
To this end, we divide the proof into several steps.

In Subsection~\ref{sect;modified-Lyapunov}, we introduce a suitably reconstructed
Lyapunov functional and derive regularity estimates for the associated
dissipative terms.
In Subsection~\ref{sect;exterior}, we establish various regularity estimates
in exterior domains by using the estimates obtained in
Subsection~\ref{sect;modified-Lyapunov}.
Thanks to the results in the previous subsections, we then derive interior
regularity estimates in Subsection~\ref{sect;interior}.
Finally, in Subsection~\ref{sect;main}, we complete the proof of
Theorem~\ref{thm;global} by combining the exterior and interior estimates.

Throughout this section, we denote by $(u,v)$ the solution to~\eqref{eqn;KS}
with maximal existence time $T\le\infty$.
In the following, we derive a priori estimates on $[t_0,\tau]\cap(0,T)$ for arbitrary $t_0,\tau$ with $t_0<\tau$, which will allow us to conclude global existence by a contradiction argument.
Moreover, $C$ denotes a generic positive constant depending only on $\lambda$
and the initial data~$(u_0,v_0)$, which may vary from line to line.
The dependence of $C$ on additional parameters will be indicated explicitly.

\subsection{Refined Lyapunov  functional}\label{sect;modified-Lyapunov}

We here introduce a suitably reconstructed Lyapunov functional which is non-increasing along the evolution of solutions, and derive regularity estimates associated with the corresponding dissipation mechanism. These estimates further yield refined regularity properties for solutions to \eqref{eqn;KS}.

The existence of a Lyapunov functional is known to play a fundamental role in the analysis of the long-time behavior of solutions, including global existence and finite-time blowup phenomena. In the whole space setting, however, the use of Lyapunov-based methods in~\eqref{eqn;KS} faces intrinsic difficulties due to 
 the necessity of controlling the behavior as $|x|\to\infty$. This issue becomes particularly delicate at the critical mass.

First of all, it is worth mentioning that the system \eqref{eqn;KS} has a usual Lyapunov functional~$\L(t)$ as follows:
\begin{align*}
\frac{\d}{\d t}\L(t)+ \D(t) =0,
\eqntag
\label{eqn;Lyapunov}
\end{align*}
where
\begin{align*}
\L(t):=\int_{\R^2}u\ln u\dx-\int_{\R^2}uv\dx+\frac12\|\nabla v\|_2^2+\frac\lambda2\|v\|_2^2
\end{align*}
and the dissipative term $\D(t)$ is given by
\begin{align*}
\D(t):=\int_{\R^2} u \left|\nabla(\ln u -v)\right|^2\dx+\|\partial_t v\|_2^2.
\end{align*}
The upper bound for $\L(t)$ follows from the monotonicity \eqref{eqn;Lyapunov} and 
the lower bound for that can be obtained by the optimal constraint of the initial mass $\|u_0\|_1$ deriving from the Trudinger--Moser type inequality and its best possible constant, see \cite{Na-Se-Yo} and also \cite{CaCo08,Mi13,NaOg11} for instance.
On the other hand,
in order to ensure the well-definedness of a usual entropy $\int_{\R^2}u\ln u \dx$ on the whole space,
at least the logarithmic moment assumption for the initial data $u_0\ln(1+|x|^2)\in L^1(\R^2)$ is required to control the behavior of solutions as $|x|\to\infty$, that is,
\begin{align*}
-\int_{\R^2}u(\ln u)_-\dx \le\,C\int_{\R^2}u\ln(1+|x|^2)\dx,
\end{align*}
where $(\ln u)_-:=\min\{\ln u,0\}$, according to \cite[Lemma 2.4]{CaCo08}.
Later, to get rid of any moment assumptions,
the modified entropy is adapted by Nagai \cite{Na11};
\begin{align*}
\int_{\R^2}(1+u)\ln(1+u)\dx,
\end{align*}
as well as the following identity for the modified Lyapunov functional $\L_m(t)$ is introduced:
\begin{align*}
\frac{\d}{\d t} \L_m(t) + \D_m(t)=\frac14\|\nabla v(t)\|_2^2,
\eqntag
\label{eqn;usualmodifiedfunctiona}
\end{align*}
where
\begin{align*}
\L_m(t):=\int_{\R^2}(1+u)\ln(1+u)\dx-\int_{\R^2}uv\dx+\frac12\|\nabla v\|_2^2+\frac\lambda2\|v\|_2^2
\end{align*}
and
\begin{align*}
\D_m(t):=\int_{\R^2}u\left|\nabla(\ln(1+u)-v)\right|^2\dx+\int_{\R^2}\left|\nabla\left(\ln(1+u)-\frac12v\right)\right|^2\dx+\|\partial_tv\|_2^2.
\end{align*}
The functional $\L_m(t)$ is no longer the Lyapunov functional, however, it is still useful for proving 
the global existence of the solution to \eqref{eqn;KS} with sub-critical mass. Indeed, it follows
from the lower bound for $\L_m(t)$ that
for any $\al\in[0,1)$
\begin{align*}
\al\int_{\R^2}(1+u)\ln(1+u)\dx+\frac12\left[1-\frac{\|u_0\|_1}{8\pi(1-\al)}\right]&\|\nabla v\|_2^2
+\int_0^t \D_m(s)\,\d s
\\
\le\,&\L_m(0)+C(\tau,\alpha)+\frac14\int_0^t\|\nabla v\|_2^2\,\d s
\eqntag
\label{eqn;afford}
\end{align*}
for $t\in[0,\tau]\cap(0,T)$ with any $\tau>0$ and some $C(\tau,\alpha)>0$ (cf.~\cite{CaCo08,Mi13}). Hence, if $\|u_0\|<8\pi$, then all terms on the left hand side of \eqref{eqn;afford} are positive by choosing $\al>0$ sufficiently small which depends on $\|u_0\|_1$, and Gronwall's inequality implies that
$\nabla v \in L^2((0,T)\times \R^2)$, 
so that we also obtain the bound for the modified entropy $\int_{\R^2}(1+u)\ln (1+u)\dx$. 
Nevertheless, for the critical mass $\|u_0\|_1=8\pi$, this necessarily leads to $\alpha=0$, as a result, 
\begin{align*}
\int_0^t \D_m(s)\,\d s\le\,&\L_m(0)+C(\tau)+\frac14\int_0^t\|\nabla v\|_2^2\,\d s,
\eqntag
\label{eqn;lack-lyapunov}
\end{align*}
which is useless and fails to obtain even the estimates for $\D_m(t)$ due to the lack of the regularity of solutions.
Therefore, we introduce a reconstructed Lyapunov functional $\F_m(t)$ so as to show regularity estimates corresponding to the dissipative terms. Let $\F_m(t)$ be the functional for solutions to \eqref{eqn;KS} defined as
\begin{align*}
\F_m(t):=\L_m(t)
+\int_{\R^2}\ln(1+u)\,\mathrm{d}x-\int_{\R^2}v\,\d x.
\eqntag
\label{eqn;modifiedLyapunov}
\end{align*}
Then, the following functional differential inequality holds true:
\begin{prop}\label{prop;energy-est}
Let $(u,v)$ be the solution to \eqref{eqn;KS}. Then, the functional $\F_m(t)$ defined in~\eqref{eqn;modifiedLyapunov} satisfies the following identity:
	\begin{align*}
	\frac{\d}{\d t}\F_m(t) +\int_{\R^2}u\left|\nabla\left(\ln(1+u)-v\right)\right|^2\,\mathrm{d}x+\frac12\|\partial_tv\|_2^2
	\le\,&
	\lambda\int_{\R^2}v\,\d x.
	\end{align*}
	Therefore, if $\lambda=0$ then
	\begin{align*}
	\frac{\d}{\d t}\F_m(t) +\int_{\R^2}u\left|\nabla\left(\ln(1+u)-v\right)\right|^2\,\mathrm{d}x+\frac12\|\partial_tv\|_2^2 \le\,0.
	\end{align*}
	If $\lambda>0$ then 
	\begin{align*}
	\frac{\d}{\d t}\F_m(t) +\int_{\R^2}u\left|\nabla\left(\ln(1+u)-v\right)\right|^2\,\mathrm{d}x+\frac12\|\partial_tv\|_2^2 \le\,\|v_0\|_1+\|u_0\|_1.
	\end{align*}
\end{prop}
\vspace{5mm}
It is worth emphasizing that $\F_m(t)$ is non-increasing in time when $\lambda=0$ as well as
unlike the identity \eqref{eqn;usualmodifiedfunctiona} based on the usual modified functional $\mathcal{L}_m(t)$, the error term appearing on the right-hand side of Proposition~\ref{prop;energy-est} can be easily controlled by initial data.
This allows us to show the regularity estimates for the dissipative term associated with $\F_m(t)$, see Proposition~\ref{prop;bound-energy} below.
\begin{lem}\label{lem;ModifiedLF}
	Suppose assumptions as in Theorem \ref{thm;global}. Then,
\begin{align*}
\frac{\d}{\d t}\F_m(t)+ \widetilde{\D}(t)=-\int_{\R^2}\partial_t v \frac{u}{1+u}\,\d x-\int_{\R^2}\frac{u}{1+u}\,\mathrm{d}x+\lambda\int_{\R^2}\frac{v}{1+u}\,\d x,
\end{align*}
where
$\F_m(t)$ is the modified functional defined in \eqref{eqn;modifiedLyapunov} and the dissipative term $\widetilde{\D}(t)$ is given by
\begin{align*}
\widetilde{\D}(t):=\int_{\R^2}u\left|\nabla \left(\ln(1+u)-v\right)\right|^2\,\d x+\|\partial_t v\|_2^2.
\end{align*}
\end{lem}

\begin{pr}{Lemma \ref{lem;ModifiedLF}}
The following computations are already well-known:
\begin{align*}
\frac{\d}{\d t} \int_{\R^2}(1+u)\ln (1+u)\,\d x=\,&
-\int_{\R^2}(1+u)|\nabla\ln(1+u)|^2\,\d x
+\int_{\R^2}u\nabla\ln(1+u)\cdot\nabla v\,\d x
\end{align*}
as well as
\begin{align*}
-\frac{\d}{\d t}\int_{\R^2}uv\,\d x=\,&\int_{\R^2}(1+u)\nabla\ln(1+u)\cdot\nabla v\,\d x-\int_{\R^2}u|\nabla v|^2\,\d x
\\
&-\|\partial_tv\|_2^2-\frac12\frac{\d}{\d t}\|\nabla v\|_2^2-\frac\lambda2\frac{\d}{\d t}\|v\|_2^2,
\end{align*}
so that
\begin{align*}
&\frac{\d}{\d t}\left[\int_{\R^2}(1+u)\ln(1+u)\,\d x-\int_{\R^2}uv\,\d x+\frac12\left(\|\nabla v\|_2^2+\lambda\|v\|_2^2\right)\right]
\\
&+\int_{\R^2}u\left|\nabla\left(\ln(1+u)-v\right)\right|^2\,\mathrm{d}x+\|\partial_tv\|_2^2
\\
=\,&-\|\nabla\ln(1+u)\|_2^2+\int_{\R^2}\nabla\ln(1+u)\cdot\nabla v\,\d x.
\end{align*}
Next,
\begin{align*}
\frac{\d}{\d t}\int_{\R^2}\ln(1+u)\,\d x=\,&\int_{\R^2}\frac{1}{1+u} \left[\Delta u -\nabla\cdot( u\nabla v)\right]\,\mathrm{d}x
\\
=\,&-\int_{\R^2}\nabla\left(\frac{1}{1+u}\right)\cdot \nabla u\,\d x+\int_{\R^2}u\nabla\left(\frac{1}{1+u}\right)\cdot \nabla v\,\d x
\\
=\,&\int_{\R^2}|\nabla \ln(1+u)|^2\,\d x-\int_{\R^2}\frac{u}{1+u}\nabla\ln(1+u)\cdot\nabla v\,\d x
\\
=\,&\int_{\R^2}|\nabla \ln(1+u)|^2\,\d x-\int_{\R^2}\nabla\ln(1+u)\cdot\nabla v\,\d x
\\
&+\int_{\R^2}\frac{1}{1+u}\nabla\ln(1+u)\cdot\nabla v \,\d x.
\end{align*}
Hence, combining the above computations implies that
\begin{align*}
&\frac{\d}{\d t}\left[\int_{\R^2}(1+u)\ln(1+u)\,\d x-\int_{\R^2}uv\,\d x+\frac12\left(\|\nabla v\|_2^2+\lambda\|v\|_2^2\right)+\int_{\R^2}\ln(1+u)\,\d x\right]
\\
&+\int_{\R^2}u\left|\nabla\left(\ln(1+u)-v\right)\right|^2\,\mathrm{d}x+\|\partial_tv\|_2^2
\\
=\,&\int_{\R^2}\frac{1}{(1+u)^2}\nabla u\cdot\nabla v \,\d x.
\end{align*}
Since
\begin{align*}
\int_{\R^2}\frac{1}{(1+u)^2}\nabla u\cdot\nabla v \,\d x=\,&-\int_{\R^2}\nabla\left(\frac{1}{1+u}\right)\cdot \nabla v\,\d x
\\
=\,&\int_{\R^2}\frac{1}{1+u}\Delta v\,\d x
\\
=\,&\int_{\R^2}\frac{\partial_t v}{1+u}\,\d x-\int_{\R^2}\frac{u}{1+u}\,\d x+\lambda\int_{\R^2}\frac{v}{1+u}\,\d x
\\
=\,&\frac{\d}{\d t}\int_{\R^2}v\,\d x-\int_{\R^2}\partial_t v \frac{u}{1+u}\,\d x-\int_{\R^2}\frac{u}{1+u}\,\d x+\lambda\int_{\R^2}\frac{v}{1+u}\,\d x,
\end{align*}
we end up with
\begin{align*}
&\frac{\d}{\d t}\F_m(t)+\int_{\R^2}u\left|\nabla\left(\ln(1+u)-v\right)\right|^2\,\mathrm{d}x+\|\partial_tv\|_2^2
\\
=\,&-\int_{\R^2}\partial_t v \frac{u}{1+u}\,\d x-\int_{\R^2}\frac{u}{1+u}\,\d x+\lambda\int_{\R^2}\frac{v}{1+u}\,\d x,
\end{align*}
as desired.
\end{pr}

\vspace{5mm}
We are now in a position to show the proof of Proposition \ref{prop;energy-est}.
\vspace{5mm}

\begin{pr}{Proposition \ref{prop;energy-est}}
By Young's inequality, 
\begin{align*}
\left|\int_{\R^2}\partial_t v \frac{u}{1+u}\,\d x \right|\le\,&\frac12\|\partial_t v\|_2^2+\frac{1}{2}\int_{\R^2}\frac{u^2}{(1+u)^2}\,\d x
\\
\le\,&\frac12\|\partial_t v\|_2^2+\frac{1}{2}\int_{\R^2}\frac{u}{1+u}\,\d x,
\end{align*}
where we use $x/(1+x)\le 1$ for $x\ge0$.
In addition, when $\lambda>0$
\begin{align*}
\lambda\int_{\R^2}\frac{v}{1+u}\,\d x\le\,\lambda \int_{\R^2}v\,\d x=\,& \lambda e^{-\lambda t} \int_{\R^2}v_0\,\d x+(1-e^{-\lambda t})\|u_0\|_1
\\
\le\,&\|v_0\|_1+\|u_0\|_1.
\end{align*}
Therefore, Lemma \ref{lem;ModifiedLF} implies the conclusion.
\end{pr}

\vspace{5mm}
By use of Proposition \ref{prop;energy-est} with $\|u_0\|_1=8\pi$,
the following regularity estimate is investigated.
\begin{prop}\label{prop;bound-energy}
Suppose that $\|u_0\|_1= 8\pi$.
Then, for any $t_0,\tau>0$ with $t_0<\tau$, there exists a constant $C(t_0, \tau)=C(t_0,\tau,\|u_0\|_1,\|v_0\|_1)>0$ such that
for any $t\in[t_0,\tau]\cap(0,T)$
\begin{align*}
\F_m(t) \ge -C(t_0,\tau).
\end{align*}
This implies that 
\begin{align*}
\int_{t_0}^t \int_{\R^2} u\left|\nabla\left(\ln (1+u)-v\right)\right|^2\,\d x\d s+\int_{t_0}^t \|\partial_t v\|_2^2\,\d s\le\, \F_m(t_0)+C(t_0,\tau),\quad t\in[t_0,\tau]\cap(0,T).
\end{align*}
\end{prop}
\vspace{5mm}
As mentioned above, such an estimate in Proposition~\ref{prop;bound-energy} cannot be expected from the usual modified functional~$\mathcal{L}_m(t)$ defined in~\eqref{eqn;usualmodifiedfunctiona}, due to the lack of regularity of solutions,  see~\eqref{eqn;lack-lyapunov}.
In the case of the critical mass, we also do not expect to obtain an estimate for the modified entropy of $u$. Nevertheless,
the functional $\F_m(t)$ introduced enjoys that estimates for such dissipative terms play an important role in establishing further a priori estimates in the subsequent analysis, see Subsection~\ref{sect;exterior} below.

In order to prove Proposition \ref{prop;bound-energy}, let us recall a functional inequality deriving from the Trudinger--Moser inequality, see for instance~\cite{Na-Se-Yo}:
\begin{lem}\label{lem;TM_dilechlet}
	Let $D$ be a domain in $\R^2$ with $|D|<\infty$.
	Then, there exists a constant $C_{TM}>0$ independent of $D$ such that for $f\in H^1_0(D)$,
	\begin{align*}
	\int_{D}\exp\left(|f|\right)\,\mathrm{d}x\le\,& C_{TM}|D| \exp\left(\frac{1}{16\pi}\|\nabla f\|_{L^2(D)}^2\right),
	\end{align*}
	where $16\pi$ is the best possible constant.
\end{lem}

\begin{pr}{Proposition \ref{prop;bound-energy}}
For $\alpha>0$,
\begin{align*}
\int_{\R^2}uv \ \mathrm{d}x=\,&\int_{ \{v> \alpha\} }uv \ \mathrm{d}x+\int_{ \{v\le \alpha\} }uv\ \mathrm{d}x
\\
\le\,&\int_{ \{v> \alpha\} }u(v-\alpha)_+ \ \mathrm{d}x+ \alpha\int_{\R^2} u \ \mathrm{d}x,
\end{align*}
where $(v-\alpha)_+:=\max\{(v-\alpha),0\}$.
Following the similar argument to \cite[Lemma 2.1]{Mi13}, we have
\begin{align*}
&\int_{ \{v\ge \alpha\} }u(v-\alpha)_+ \ \mathrm{d}x
\\
\le\,& \int_{ \{v> \alpha\} }u\ln u\ \mathrm{d}x
+ \left(\int_{ \{v> \alpha\} }u\,\d x\right)\ln\left(\int_{ \{v> \alpha\} }\exp\left((v-\alpha)_+\right)\ \mathrm{d}x\right)
\\
&-\left(\int_{ \{v> \alpha\} }u\,\d x\right)\ln \left(\int_{ \{v> \alpha\} }u\,\d x\right)
\\
\le\,& \int_{\R^2 }(1+u)\ln (1+u)\ \mathrm{d}x
+ \|u_0\|_1\ln\left(\int_{ \{v>\alpha\} }\exp\left((v-\alpha)_+\right)\ \mathrm{d}x\right)+\frac1e,
\end{align*}
where we use 
\begin{align*}
\int_{D} gh\ \mathrm{d}x\le \int_{D} g\log g \ \mathrm{d}x+ M \log \left(\int_{D}e^h \ \mathrm{d}x\right)-M\log M, 
\quad M:=\int_{D} g\ \mathrm{d}x,
\eqntag
\label{eqn;NSS}
\end{align*}
see \cite[Lemma 5.4]{NaSeSu} and \cite[Proposition 2.3]{Mi13}.
Since $(v-\alpha)_+ \in H^1_0(\{v>\alpha\})$, it follows from the two-dimensional Trudinger--Moser type inequality stated in Lemma \ref{lem;TM_dilechlet} that
\begin{align*}
\int_{ \{v>\alpha\} }\exp\left((v-\alpha)_+\right)\ \mathrm{d}x
\le\,& C_{TM}|\{v>\alpha\}| \exp\left(\frac{\|\N v\|_2^2}{16\pi } \right),
\end{align*}
which implies that
\begin{align*}
\|u_0\|_1\ln\left(\int_{ \{v>\alpha\} }\exp\left((v-\alpha)_+\right)\ \mathrm{d}x\right)
\le \frac{\|u_0\|_1}{16\pi } \|\N v\|_2^2+ \|u_0\|_1 \log (C_{TM}|\{v>\alpha\}|),
\end{align*}
so that
\begin{align*}
\int_{\R^2 }uv\,\d x\le\,& \int_{\R^2 }(1+u)\ln (1+u)\ \mathrm{d}x
\\&+\frac{\|u_0\|_1}{16\pi } \|\N v\|_2^2+ \|u_0\|_1 \log (C_{TM}|\{v>\alpha\}|)+\frac1e+\alpha\|u_0\|_1
\eqntag
\label{eqn;L^1-uv}
\end{align*}
for any $\alpha>0$.
We thus obtain from \eqref{eqn;L^1-uv} along with $\|u_0\|_1=8\pi$
\begin{align*}
\F_m(t):=\,&\int_{\R^2}(1+u)\ln(1+u)\,\d x-\int_{\R^2}uv\,\d x+\frac12\left(\|\nabla v\|_2^2+\lambda\|v\|_2^2\right)
\\&+\int_{\R^2}\ln(1+u)\,\d x-\int_{\R^2}v\,\d x
\\
\ge\,&\frac\lambda2\|v\|_2^2+\int_{\R^2}\ln(1+u)\,\d x
-\int_{\R^2}v\,\d x- \|u_0\|_1 \ln (C_{TM}|\{v>\alpha\}|)-\frac1e-\alpha\|u_0\|_1
\\
\ge\,&-\int_{\R^2}v\,\d x- \|u_0\|_1 \ln \left(\frac{C_{TM}\|v\|_1}{\alpha}\right)-\frac1e-\alpha\|u_0\|_1.
\end{align*}
Finally, since it follows from the second equation to \eqref{eqn;KS} that
\begin{align*}
\int_{\R^2 }v\,\d x
=
\left\{
\begin{aligned}
&\|v_0\|_1+ t \|u_0\|_1 &\text{if}~\lambda=0,
\\
&e^{-\lambda t}\|v_0\|_1+\frac{1-e^{-\lambda t}}{\lambda}\|u_0\|_1&\text{if}~\lambda>0,
\end{aligned}
\right.
\end{align*}
we end up with the lower bound for $\F_m(t)$ that
for any $t_0>0$ and $\tau>0$ with $t_0<\tau$, there exists a constant $C(t_0,\tau)>0$ such that
\begin{align*}
\F_m(t) \ge\, -C(t_0,\tau)
\end{align*}
for $t\in[t_0,\tau]\cap(0,T)$.
Next, by Proposition \ref{prop;energy-est}
\begin{align*}
\F_m(t)+
\int_{t_0}^t \int_{\R^2} u\left|\nabla\left(\ln (1+u)-v\right)\right|^2\,\d x\d s+\frac12\int_{t_0}^t \|\partial_t v\|_2^2\,\d s
\le\,&\F_m(t_0)+\lambda \int_{t_0}^t\int_{\R^2}v\,\mathrm{d}x\d s
\end{align*}
for $t\in[t_0,\tau]\cap(0,T)$.
Using the lower bound for $\F_m(t)$, we derive  the desired estimate that
\begin{align*}
\int_{t_0}^t \int_{\R^2} u\left|\nabla\left(\ln (1+u)-v\right)\right|^2\,\d x\d s+\frac12\int_{t_0}^t \|\partial_t v\|_2^2\,\d s
\le\,&\F_m(t_0)+ C(t_0,\tau),
\end{align*}
which concludes the proof.
\end{pr}

\vspace{5mm}
\subsection{Regularity estimates in  exterior domains}\label{sect;exterior}

This subsection is devoted to the study of regularity estimates for solutions in exterior domains.
In the critical mass case for the fully parabolic system \eqref{eqn;KS} in $\R^2$, we  realize that the bounds on the dissipative terms derived from the refined Lyapunov functional $\F_m(t)$ introduced in Subsection~\ref{sect;modified-Lyapunov} play an essential role in establishing further regularity estimates in exterior domains, which in turn lead to the global existence of solutions to \eqref{eqn;KS}.

Define the cut-off function $\varphi$  in $C^\infty([0,\infty))$ by
\begin{align*}
\varphi=\varphi(|x|)
=
\left\{
\begin{aligned}
&0,&&|x|<1
\\
&\text{smooth},&&1\le|x|\le 2,
\\
&1,&&|x|>2,
\end{aligned}
\right.
\eqntag
\label{eqn;cut-off}
\end{align*}
with $|\nabla \varphi|, |\nabla^2\varphi|\le\, C\varphi$ for some constant $C>0$.
Set $\varphi_R=\varphi_R(|x|)=\varphi(|x|/R)$ for $R>0$.

\begin{prop}\label{prop;exterior-L^1}
Suppose that $\|u_0\|_1=8\pi$.
For each $t_0,\tau>0$ with $t_0<\tau$, there is a constant $C(t_0,\tau)>0$ such that
for every $R>0$
\begin{align*}
\int_{\R^2}u\varphi_{\frac R2}\,\d x\le\,&\int_{\R^2}u(t_0)\varphi_{\frac R2}\,\d x+\frac{C(t_0,\tau)}{R}+\frac{C(t_0,\tau)}{R^2}
\end{align*}
for $t\in[t_0,\tau]\cap(0,T)$. This consequence implies that
\begin{align*}
\int_{ |x|>R } u\dx\le\,&\int_{\R^2}u(t_0)\varphi_{\frac R2}\,\d x+\frac{C(t_0,\tau)}{R}+\frac{C(t_0,\tau)}{R^2}.
\end{align*}
\end{prop}

\begin{rem}\label{rem;L^1-exterior}
	Proposition \ref{prop;exterior-L^1} gives the fact that for any $\ep>0$, there is $R_0=R_0(\ep,t_0,\tau)>0$ sufficiently large such that
	for any $R>R_0$
	\begin{align*}
\int_{|x|>R}u\,\d x <\ep
	\end{align*}
	for $t\in[t_0,\tau]\cap(0,T)$
	since it  follows by the dominated convergence theorem that
	\begin{align*}
	\lim_{R\to\infty}\int_{\R^2}u(t_0)\varphi_{\frac R2}\,\d x=0.
	\end{align*}
\end{rem}

\begin{rem}
According to \cite[Eq.~(3.21)]{NaOg16} and \cite[Lemma~3.2]{NaYa20NA}, the counterpart of Proposition~\ref{prop;exterior-L^1} for the parabolic-elliptic system, which is the system \eqref{eqn;PKS}, follows readily from the well-known integral symmetry based on the elliptic operator $(-\Delta+\lambda)^{-1}$ acting on the second component.
By contrast, establishing the analogous result for the fully parabolic Keller--Segel system \eqref{eqn;KS} in $\R^2$ is a nontrivial task under the critical mass condition $\|u_0\|_1=8\pi$, particularly in the absence of any moment assumptions on the initial data, due to the lack of regularity of solutions.
\end{rem}

\begin{pr}{Proposition \ref{prop;exterior-L^1}}
Differentiating the $L^1$-norm of $u\varphi_{R/2}$ with respect to $t$ and integrating by parts, we have
\begin{align*}
\frac{\d}{\d t}\int_{\R^2}u\varphi_{\frac R2}\,\d x=\,&-\int_{\R^2}u\nabla\left(\ln u -v\right)\cdot \nabla\varphi_{\frac R2}\,\d x
\\
=\,&-\int_{\R^2}u\nabla\left(\ln (1+u) -v\right)\cdot \nabla\varphi_{\frac R2}\,\d x
-\int_{\R^2}u\nabla\ln u \cdot \nabla \varphi_{\frac R2}\,\d x
\\
&+\int_{\R^2}u\nabla\ln (1+u) \cdot \nabla \varphi_{\frac R2}\,\d x
\\
=\,&-\int_{\R^2}u\nabla\left(\ln (1+u) -v\right)\cdot \nabla\varphi_{\frac R2}\,\d x-\int_{\R^2}\nabla\ln (1+u) \cdot \nabla \varphi_{\frac R2}\,\d x
\\
=\,&-\int_{\R^2}u\nabla\left(\ln (1+u) -v\right)\cdot \nabla\varphi_{\frac R2}\,\d x
+\int_{\R^2}\ln(1+u) \Delta \varphi_{\frac R2}\,\d x.
\end{align*}
Notice that $|\nabla\varphi_{R/2}|\le\,C\varphi_{R/2}/R$ and $|\Delta\varphi_R|\le\,C\varphi_{ R/2}/R^2$, so that
Young's inequality along with the mass conservation implies that
\begin{align*}
\left| \int_{\R^2}u\nabla\left(\ln (1+u) -v\right)\cdot \nabla\varphi_{\frac R2}\,\d x \right|
\le\,&\frac{C}{R}\int_{\R^2}u|\nabla(\ln (1+u)-v)| \varphi_{\frac R2}\,\d x
\\
\le\,&\frac{C}{R}\int_{\R^2}u\varphi_{\frac R2}\,\d x+\frac{C}{R}\int_{\R^2}u|\nabla(\ln(1+u)-v)|^2\varphi_{\frac R2}\,\d x
\\
\le\,&\frac{C\|u_0\|_1}{R}+\frac{C}{R}\int_{\R^2}u|\nabla(\ln(1+u)-v)|^2\,\d x
\end{align*}
as well as
\begin{align*}
\left| \int_{\R^2}\ln(1+u) \Delta \varphi_{\frac R2}\,\d x \right|
\le\,&\frac{C}{R^2}\int_{\R^2}\ln(1+u)\varphi_{\frac R2}\,\d x
\le\, \frac{C\|u_0\|_1}{R^2}.
\end{align*}
Hence, integrating this over $(t_0,t)$ leads us to
\begin{align*}
\int_{\R^2}u\varphi_{\frac R2}\,\d x\le\,&\int_{\R^2}u(t_0)\varphi_{\frac R2}\,\d x+
\frac{C\|u_0\|_1}{R}t+\frac{C\|u_0\|_1}{R^2}t+\frac{C}{R}\int_{t_0}^t\int_{\R^2}u|\nabla(\ln(1+u)-v)|^2\,\d x\d s
\\
\le\,&\int_{\R^2}u(t_0)\varphi_{\frac R2}\,\d x+
\frac{C\|u_0\|_1}{R}\tau+\frac{C\|u_0\|_1}{R^2}\tau+\frac{C}{R}\int_{t_0}^t\int_{\R^2}u|\nabla(\ln(1+u)-v)|^2\,\d x\d s
\end{align*}
for $t\in[t_0,\tau]\cap(0,T)$ with any $t_0,\tau>0$.
By virtue of Proposition \ref{prop;bound-energy},
\begin{align*}
\int_{t_0}^t\int_{\R^2}u|\nabla(\ln(1+u)-v)|^2\,\d x\d s\le\,& C(t_0,\tau).
\end{align*}
Therefore, we have
\begin{align*}
\int_{\R^2}u\varphi_{\frac R2}\,\d x\le\,&\int_{\R^2}u(t_0)\varphi_{\frac R2}\,\d x+
\frac{C\|u_0\|_1}{R}t+\frac{C\|u_0\|_1}{R^2}t+\frac{C}{R}\int_{t_0}^t\int_{\R^2}u|\nabla(\ln(1+u)-v)|^2\,\d x\d s
\\
\le\,&\int_{\R^2}u(t_0)\varphi_{\frac R2}\,\d x+
\frac{C(t_0,\tau)}{R}+\frac{C(t_0,\tau)}{R^2}
\end{align*}
for $t\in[t_0,\tau]\cap(0,T)$. 
In addition, since $\varphi_{ R/2}\equiv1$ on $|x|\ge R$,
\begin{align*}
\int_{ |x|>R } u\dx \le\,\int_{\R^2 }u\varphi_{\frac R2}\dx.
\end{align*}
The proof is complete.
\end{pr}

\vspace{5mm}
Thanks to Lemma \ref{prop;exterior-L^1}, the constraint mass on the exterior domain can be taken arbitrary small as desired depending on $R>0$. 
This leads us to
 the further regularity estimates.
\begin{lem}\label{lem;LlogL-bound-exterior}
Suppose assumptions as in Theorem \ref{thm;global}.
For any $t_0,\tau>0$ with $t_0<\tau$, there is a constant $R_0=R_0(t_0,\tau,\|u_0\|_1)\gg1$ sufficiently large such that
for all $R>R_0$,
there is a constant $C(t_0,\tau,R)>0$ satisfying
\begin{align*}
\int_{ |x|>2R }(1+u)\ln(1+u)\dx +\int_{t_0}^t \int_{\R^2}|\nabla v|^2\varphi_R\dx\d s\le\, C(t_0,\tau,R)
\end{align*}
for $t\in[t_0,\tau]\cap(0,T)$.
\end{lem}

\begin{pr}{Lemma \ref{lem;LlogL-bound-exterior}}
Let $\varphi_R$ be the cut-off function defined in \eqref{eqn;cut-off}.
By integrations by parts, we have
\begin{align*}
&\frac{\d}{\d t}\left(\int_{\R^2}(1+u)\ln(1+u)\varphi_R\dx-\int_{\R^2}u\varphi_R\dx\right)
\\
=\,&\int_{\R^2}\ln(1+u)\varphi_R\left[\Delta u - \nabla\cdot(u\nabla v)\right]\,\dx
\\
=\,&-\int_{\R^2 }\frac{|\nabla u|^2}{1+u}\varphi_R\dx-\int_{\R^2}\ln(1+u)\nabla u \cdot \nabla\varphi_R\dx
\\
&+\int_{\R^2 }\frac{u}{1+u}\varphi_R\nabla u\cdot\nabla v \dx+\int_{\R^2}u\ln(1+u)\nabla v \cdot \nabla\varphi_R\dx.
\end{align*}
In light of the second equation in \eqref{eqn;KS},
\begin{align*}
\frac12\frac{\d}{\d t}\int_{\R^2}v^2\varphi_R\dx=\,&\int_{\R^2 }v\varphi_R\left[\Delta v -\lambda v +u\right]\dx
\\
=\,&-\int_{\R^2}|\nabla v|^2\varphi_R\dx-\int_{\R^2}v \nabla v \cdot\nabla  \varphi_R\dx-\lambda\int_{\R^2}v^2\varphi_R\dx+\int_{\R^2 }uv\varphi_R\dx,
\end{align*}
so that
\begin{align*}
&\frac{\d}{\d t}\left(\int_{\R^2}(1+u)\ln(1+u)\varphi_R\dx-\int_{\R^2}u\varphi_R\dx+\frac12\int_{\R^2 }v^2\varphi_R\dx\right)
\\
=\,&-\int_{\R^2}\frac{|\nabla u|^2}{1+u}\varphi_R\dx-\int_{\R^2}|\nabla v|^2\varphi_R\dx
-\int_{\R^2}\ln(1+u)\nabla u \cdot \nabla\varphi_R\dx
\\
&+\int_{\R^2 }\frac{u}{1+u}\varphi_R\nabla u\cdot\nabla v \dx+\int_{\R^2}u\ln(1+u)\nabla v \cdot \nabla\varphi_R\dx
\\
&-\int_{\R^2}v \nabla v \cdot \nabla \varphi_R\dx-\lambda\int_{\R^2}v^2\varphi_R\dx+\int_{\R^2 }uv\varphi_R\dx.
\eqntag
\label{eqn;LlogL1}
\end{align*}
By Young's inequality and $\ln(1+x)\le\, x/(\sqrt{1+x})$ for $x\ge0$, we see for any $\ep>0$
\begin{align*}
\left| \int_{\R^2}\ln(1+u)\nabla u \cdot \nabla\varphi_R\dx \right|
\le\,&\ep\int_{\R^2}\frac{|\nabla u|^2}{1+u}\varphi_R\dx+\frac{C_\ep}{R^2}\int_{\R^2 }(1+u)\left(\ln(1+u)\right)^2\varphi_R\dx
\\
\le\,&\ep\int_{\R^2}\frac{|\nabla u|^2}{1+u}\varphi_R\dx+\frac{C_\ep}{R^2}\int_{\R^2 }u^2\varphi_R\dx
\end{align*}
and
\begin{align*}
&\int_{\R^2}u\ln(1+u)\nabla v \cdot \nabla\varphi_R\dx
\\
=\,&-\int_{\R^2}u\ln(1+u)\nabla(\ln(1+u)-v)\cdot\nabla \varphi_R\dx
+\int_{\R^2}\frac{u}{1+u}\ln(1+u) \nabla u \cdot \nabla \varphi_R\dx
\\
\le\,&\frac{C}{R^2}\int_{\R^2 } u(\ln(1+u))^2\varphi_R\dx+\int_{\R^2}u\left|\nabla(\ln(1+u)-v)\right|^2\varphi_R\dx
\\
&+\ep\int_{\R^2}\frac{|\nabla u|^2}{1+u}\varphi_R\dx+\frac{C_\ep}{R^2}\int_{\R^2 }\frac{u^2}{1+u}\left(\ln(1+u)\right)^2\varphi_R\dx
\\
\le\,&\ep\int_{\R^2}\frac{|\nabla u|^2}{1+u}\varphi_R\dx
+\frac{C_\ep}{R^2}\int_{\R^2 } u^2\varphi_R\dx+\int_{\R^2}u\left|\nabla(\ln(1+u)-v)\right|^2\dx
\end{align*}
as well as by Lemma \ref{lem;v-Lpbound}
\begin{align*}
&-\int_{\R^2}v \nabla v \cdot \nabla \varphi_R\dx-\lambda\int_{\R^2}v^2\varphi_R\dx+\int_{\R^2 }uv\varphi_R\dx
\\
\le\,&\ep \int_{\R^2}|\nabla v|^2\dx+\frac12\int_{\R^2}u^2\varphi_R\dx+C_\ep\sup_{t>0}\|v(t)\|_2^2.
\end{align*}
Using the second equation in \eqref{eqn;KS}, we derive
\begin{align*}
\int_{\R^2 }\frac{u}{1+u}\varphi_R\nabla u\cdot\nabla v \dx
=\,&\int_{\R^2}\varphi_R\nabla u\cdot \nabla v\dx-\int_{\R^2 }\varphi_R\nabla \ln(1+u)\cdot \nabla v\dx
\\
=\,&-\int_{\R^2 }u\varphi_R\left[\partial_t v +\lambda v -u\right]\dx-\int_{\R^2 }u\nabla v \cdot \nabla \varphi_R\dx
\\
&-\int_{\R^2 }\varphi_R\nabla \ln(1+u)\cdot \nabla v\dx
\\
\le\,&\frac12\int_{\R^2}u^2\varphi_R\dx+\frac12\int_{\R^2}|\partial_t v|^2\dx-\lambda\int_{\R^2}uv\varphi_R\dx+\int_{\R^2}u^2\varphi_R\dx
\\
&+\frac{C_\ep}{R^2}\int_{\R^2}u^2\varphi_R\dx+\ep\int_{\R^2}|\nabla v|^2\varphi_R\dx
\\
&+\frac12\int_{\R^2}|\nabla\ln(1+u)|^2\varphi_R\dx+\frac12\int_{\R^2 }|\nabla v|^2\varphi_R\dx
\\
\le\,&\frac12\int_{\R^2}\frac{|\nabla u|^2}{1+u}\varphi_R\dx+\left(\frac32+\frac{C_\ep}{R^2}\right)\int_{\R^2}u^2\varphi_R\dx
\\
&+\left(\ep+\frac12\right)\int_{\R^2}|\nabla v|^2\varphi_R\dx+\frac12\int_{\R^2}|\partial_t v|^2\dx.
\end{align*}
Combining above outcomes implies from \eqref{eqn;LlogL1} that
\begin{align*}
&\frac{\d}{\d t}\left(\int_{\R^2}(1+u)\ln(1+u)\varphi_R\dx-\int_{\R^2}u\varphi_R\dx+\frac12\int_{\R^2 }v^2\varphi_R\dx\right)
\\
\le\,&-\left(1-2\ep -\frac12\right)\int_{\R^2}\frac{|\nabla u|^2}{1+u}\varphi_R\dx
-\left(1-2\ep-\frac12\right)\int_{\R^2}|\nabla v|^2\varphi_R\dx
\\
&+\left(2+\frac{C_\ep}{R^2}\right)\int_{\R^2 }u^2\varphi_R\dx
+\int_{\R^2}u\left|\nabla(\ln(1+u)-v)\right|^2\dx+\frac12\int_{\R^2}|\partial_t v|^2\dx+C_\ep.
\end{align*}
According to Lemma \ref{lem;L2andL3} as $|f|=u$ and $\phi=\varphi_R$, we see
\begin{align*}
\int_{\R^2}u^2\varphi_R\dx\le\,&2\left(\int_{ \{u>1\}\cap \{\supp\varphi_R\} }u\dx\right)\left(\int_{ \{u>1\} }\frac{|\nabla u|^2}{1+u} \varphi_R\dx\right)
\\
&+4\left(\int_{\R^2}|u\nabla \varphi_R^{\frac12}|\dx\right)^2+4\left(\int_{\R^2}u\varphi_R\dx\right)
\\
\le\,&2\left(\int_{|x|>R}u\dx\right)\left(\int_{\R^2 }\frac{|\nabla u|^2}{1+u} \varphi_R\dx\right)+C(\|u_0\|_1,R)
\end{align*}
since $\supp\varphi_R\subset \{x; |x|>R\}$. Hence,
\begin{align*}
&\frac{\d}{\d t}\left(\int_{\R^2}(1+u)\ln(1+u)\varphi_R\dx-\int_{\R^2}u\varphi_R\dx+\frac12\int_{\R^2 }v^2\varphi_R\dx\right)
\\
\le\,&-\left(\frac12-2\ep -2\left(2+\frac{C_\ep}{R^2}\right)\int_{|x|>R}u\dx\right)\int_{\R^2}\frac{|\nabla u|^2}{1+u}\varphi_R\dx-\left(\frac12-2\ep\right)\int_{\R^2}|\nabla v|^2\varphi_R\dx
\\
&
+\int_{\R^2}u\left|\nabla(\ln(1+u)-v)\right|^2\dx+\frac12\int_{\R^2}|\partial_t v|^2\dx+C_\ep
\eqntag
\label{eqn;LlogL2}
\end{align*}
for $t\in[t_0,\tau]\cap (0,T)$.
Choosing $\ep=1/8$ shows
\begin{align*}
&\frac{\d}{\d t}\left(\int_{\R^2}(1+u)\ln(1+u)\varphi_R\dx-\int_{\R^2}u\varphi_R\dx+\frac12\int_{\R^2 }v^2\varphi_R\dx\right)
\\
\le\,&-\left(\frac14-2\left(2+\frac{C}{R^2}\right)\int_{|x|>R}u\dx\right)\int_{\R^2}\frac{|\nabla u|^2}{1+u}\varphi_R\dx
-\frac14\int_{\R^2}|\nabla v|^2\varphi_R\dx
\\
&
+\int_{\R^2}u\left|\nabla(\ln(1+u)-v)\right|^2\dx+\frac12\int_{\R^2}|\partial_t v|^2\dx+C
\end{align*}
and also by virtue of Proposition \ref{prop;exterior-L^1}, one can find a constant $R_0=R_0(t_0,\tau,\|u_0\|_1)\gg1$ such that for all $R>R_0$
\begin{align*}
\frac14-2\left(2+\frac{C}{R^2}\right)\sup_{t_0\le t \le \tau}\int_{|x|>R}u(t)\dx>0.
\end{align*}
Therefore, integrating \eqref{eqn;LlogL2} over $(t_0,t)$, together with the mass conservation, Lemma~\ref{lem;v-Lpbound} and Proposition~\ref{prop;bound-energy}, leads us to
\begin{align*}
&\int_{\R^2}(1+u)\ln(1+u)\varphi_R\dx+\frac14\int_{t_0}^t \int_{\R^2}|\nabla v|^2\varphi_R\dx\d s\\
\le\,&\int_{\R^2}(1+u(t_0))\ln(1+u(t_0))\dx+\frac12\|v(t_0)\|_2^2
+\int_{t_0}^t\int_{\R^2}u\left|\nabla(\ln(1+u)-v)\right|^2\dx\d s
\\&+\int_{t_0}^t\int_{\R^2}|\partial_t v|^2\dx\d s + Ct
\\
\le\,&\int_{\R^2}(1+u(t_0))\ln(1+u(t_0))\dx+\frac12\|v(t_0)\|_2^2+C(t_0,\tau)
\end{align*}
for $t\in[t_0,\tau]\cap (0,T)$. As a result, we have
\begin{align*}
\int_{|x|>2R}(1+u)\ln(1+u)\dx\le\,\int_{\R^2}(1+u)\ln(1+u)\varphi_R\dx\le\, C(t_0,\tau)
\end{align*}
for $t\in[t_0,\tau]\cap (0,T)$. The proof is complete.
\end{pr}

\begin{lem}\label{lem;L^2-exterior}
Suppose assumptions as in Theorem \ref{thm;global}. Let $1<p <\infty$.
For any $t_0,\tau>0$ with $t_0<\tau$, let $R_0=R_0(t_0,\tau)>0$ be taken as in Lemma~\ref{lem;LlogL-bound-exterior}.
Then, there exists a constant $R_1\ge R_0>0$ such that
for all $R>R_1$, there is a constant $C(t_0,\tau,R)>0$ satisfying
\begin{align*}
\| (u \varphi_{2R})\|_2+\| (u \varphi_{4R})\|_p\le\,C(t_0,\tau,R)
\end{align*}
as well as
\begin{align*}
\int_{  |x|>4R }u^2\dx +\int_{  |x|>8R } u^p\dx \le\, C(t_0,\tau,R)
\end{align*}
for $t\in[t_0,\tau]\cap(0,T)$.
\end{lem}
\begin{pr}{Lemma \ref{lem;L^2-exterior}}
Let us compute the $L^p$-energy estimates.
By integration by parts with the first equation for \eqref{eqn;KS}, we see for $p>1$
\begin{align*}
&\frac{\d}{\d t}\int_{\R^2}(u\varphi_{2R})^p\,\d x
\\
=\,&p\int_{\R^2}u^{p-1}\varphi_{2R}^p(\Delta u -\nabla\cdot(u\nabla v))\,\d x
\\
=\,&-\frac{4(p-1)}{p}\int_{\R^2}|\nabla u^{\frac p2}|^2\varphi_{2R}^p\,\d x-2\int_{\R^2}u^{\frac p2}\nabla u^{\frac p2} \cdot\nabla\varphi_{2R}^p\dx
\\
&-(p-1)\int_{\R^2}u^p\Delta v\varphi_{2R}^p\dx+\int_{\R^2}u^p\nabla v\cdot \nabla\varphi_{2R}^p\dx
\\
=\,&-\frac{4(p-1)}{p}\int_{\R^2}|\nabla u^{\frac p2}|^2\varphi_{2R}^p\,\d x-2\int_{\R^2}u^{\frac p2}\nabla u^{\frac p2} \cdot\nabla\varphi_{2R}^p\dx
\\
&-(p-1)\int_{\R^2}u^p\partial_tv\varphi_{2R}^p\dx+(p-1)\int_{\R^2}u^{p+1}\varphi_{2R}^p\dx
-(p-1)\lambda\int_{  \R^2 }u^pv\varphi_{2R}^p\dx
\\
&-2\int_{\R^2}vu^{\frac p2}\nabla u^{\frac p2}\cdot\nabla\varphi_{2R}^p\dx
-\int_{\R^2}vu^p\Delta \varphi_{2R}^p\dx.
\eqntag
\label{eqn;L^p-exterior1}
\end{align*}
We  first derive the $L^2$-estimate for $u\varphi_{2R}$ for \eqref{eqn;L^p-exterior1} with $p=2$, that is,
\begin{align*}
\frac{\d}{\d t}\int_{\R^2}(u\varphi_{2R})^2\,\d x
=\,&-2\int_{\R^2}|\nabla u|^2\varphi_{2R}^2\,\d x-2\int_{\R^2}u\nabla u \cdot\nabla\varphi_{2R}^2\dx
\\
&-\int_{\R^2}u^2\partial_tv\varphi_{2R}^2\dx+\int_{\R^2}u^{3}\varphi_{2R}^2\dx
-\lambda\int_{  \R^2 }u^2v\varphi_{2R}^p\dx
\\
&-2\int_{\R^2}vu\nabla u\cdot\nabla\varphi_{2R}^2\dx
-\int_{\R^2}vu^2\Delta \varphi_{2R}^2\dx.
\eqntag
\label{eqn;L^2-exterior1}
\end{align*}
As for the second term on the right had side of \eqref{eqn;L^2-exterior1},
\begin{align*}
2\left| \int_{\R^2}u\nabla u \cdot \nabla \varphi_{2R}^2\dx \right|
\le\,&\frac{C}{R}\int_{\R^2}u|\nabla u| \varphi_{2R}^2\dx\\
\le\,&\ep \int_{\R^2}|\nabla u|^2\varphi_{2R}^2\dx+\frac{C_\ep}{R^2}\int_{\R^2}u^2\varphi_{2R}^2\dx
\end{align*}
for any $\ep>0$.
For the third term,
\begin{align*}
\left| \int_{\R^2}u^2 \partial_t v\varphi_{2R}^2\dx \right|
\le\,&\|u\varphi_{2R}\|_4^2\|\partial_t v\|_2
\\
\le\,&C \|u\varphi_{2R}\|_2\|\nabla(u\varphi_{2R})\|_2\|\partial_tv\|_2
\\
\le\,&C \|u\varphi_{2R}\|_2\||\nabla u|\varphi_{2R}\|_2\|\partial_tv\|_2+\frac{C}{R^2}\|u\varphi_{2R}\|_2^2\|\partial_tv\|_2
\\
\le\,&\ep\int_{\R^2}|\nabla u|^2\varphi_{2R}^2\dx+C_\ep\|u\varphi_{2R}\|_2^2\left(1+\frac{\|\partial_t v\|_2^2}{R^2}\right).
\end{align*}
The sixth term is also estimated as
\begin{align*}
2\left| \int_{\R^2}vu \nabla u\cdot \varphi_{2R}^2\dx \right|
\le\,&\frac{C}{R}\int_{\R^2}v u |\nabla u| \varphi_{2R}^2\dx
\\
\le\,&\ep\int_{\R^2}|\nabla u|^2\varphi_{2R}^2\dx+\int_{\R^2}u^3\varphi_{2R}^3\dx+\frac{C_\ep}{R^6}\int_{\R^2}v^6\dx
\end{align*}
and for the last term
\begin{align*}
\left|  \int_{\R^2}vu^2\Delta \varphi_{2R}^2\dx\right|\le\,&\frac{C}{R^2}\int_{\R^2}v u^2\varphi_{2R}^2\dx
\\
\le\,&\int_{\R^2}u^3\varphi_{2R}^3 \dx+\frac{C}{R^6}\int_{\R^2}v^3\dx,
\end{align*}
so that it follows from above outcomes that
\begin{align*}
&\frac{\d}{\d t}\int_{\R^2}(u\varphi_{2R})^2\,\d x+\left(2-3\ep\right)\int_{\R^2 }|\nabla u|^2\varphi_{2R}^2\dx
\\
\le\,&\int_{\R^2}u^3\varphi_{2R}^2\dx+2\int_{\R^2}(u\varphi_{2R})^3\dx+C_\ep\|u\varphi_{2R}\|_2^2\left(1+\frac{1+\|\partial_t v\|_2^2}{R^2}\right)+\frac{C_\ep\|v\|_6^6+C\|v\|_3^3}{R^6}.
\end{align*}
Using Lemma \ref{lem;L2andL3}, we have
\begin{align*}
\int_{\R^2}u^3\varphi_{2R}^2\dx
\le\,&\ep \left(\int_{  \{\supp\varphi_{2R}\} }(1+u)\ln(1+u)\dx\right)\left(\int_{\R^2 }|\nabla u|^2\varphi_{2R}^2\dx\right)
\\
&+C\left(\int_{\R^2}|u^{\frac32} \nabla \varphi_{2R}|\dx\right)^2+C_\ep\left(\int_{\R^2}u\varphi_{2R}^2\dx\right)
\\
\le\,&\ep \left(\sup_{t_0\le t \le \tau}\int_{  |x|>2R }(1+u(t))\ln(1+u(t))\dx\right)\left(\int_{\R^2 }|\nabla u|^2\varphi_{2R}^2\dx\right)
\\
&+\frac{C\|u_0\|_1}{R^2}\|u\varphi_{2R}\|_2^2+C_\ep\|u_0\|_1,
\end{align*}
and it follows by Gagliardo--Nirenberg's inequality that
\begin{align*}
\| u \varphi_{2R}\|_3^3\le\,&C \| u \varphi_{2R}\|_1\|\nabla(u\varphi_{2R})\|_2^2
\le\,C\|u\varphi_{2R}\|_1\||\nabla u|\varphi_{2R}\|_2^2+\frac{C\|u_0\|_1}{R^2}\|u\varphi_R\|_2^2.
\end{align*}
In addition, by the standard parabolic regularity in Lemma \ref{lem;v-Lpbound},
\begin{align*}
\|v\|_r\le C(\tau,r,\|v_0\|_r, \|u_0\|_1),\quad r\in[1,\infty).
\end{align*}
Hence,
\begin{align*}
&\frac{\d}{\d t}\int_{\R^2}(u\varphi_{2R})^2\,\d x+\left(2-3\ep-\ep C(\tau,R)-C\|u\varphi_{2R}\|_1\right)\int_{\R^2 }|\nabla u|^2\varphi_{2R}\dx
\\
\le\,&C_\ep\|u\varphi_{2R}\|_2^2\left(1+\frac{1+\|\partial_t v\|_2^2}{R^2}\right)
+\frac{C_\ep}{R^6}
\end{align*}
for $t\in[t_0,\tau]\cap(0,T)$ and any $\ep>0$.
According to Proposition \ref{prop;exterior-L^1} with Remark \ref{rem;L^1-exterior},
one can find a constant $R_1=R_1(t_0,\tau,\|u_0\|_1)\ge R_0\gg1$ sufficiently large satisfying for all $R>R_1$
\begin{align*}
2-3\ep-\ep C(t_0,\tau,R)-C\|u\varphi_{2R}\|_1>0,
\end{align*}
together with taking sufficiently small $\ep>0$.
Consequently, we have
\begin{align*}
\|u\varphi_{2R}\|_2^2\le\,&\|(u\varphi_{2R})(t_0)\|_2^2
+C\int_{t_0}^t \|u\varphi_{2R}\|_2^2\left(1+\frac{1+\|\partial_t v\|_2^2}{R^2}\right)\d s
+\frac{C}{R^6}t
\end{align*}
by integrating over $(t_0,t)$.
We already know the following bound by Proposition \ref{prop;bound-energy}
\begin{align*}
\int_{t_0}^t \|\partial_t v\|_2^2\,\d s\le\, C(t_0,\tau),\quad t\in[t_0,\tau]\cap(0,T),
\end{align*}
and we thus obtain the $L^2$-bound for $u\varphi_{2R}$ after using Gronwall's inequality that
\begin{align*}
\|u\varphi_{2R}\|_2^2\le\,&C(t_0,\tau,R) \exp\left(\int_{t_0}^t \left(1+\frac{1+\|\partial_t v\|_2^2}{R^2}\right)\d s\right)
\le\,C(t_0,\tau,R),\quad t\in[t_0,\tau]\cap(0,T).
\end{align*}
Thus, we have
\begin{align*}
\int_{  |x|>4R } u^2\dx\le\,&\int_{\R^2 }(u\varphi_{2R})^2\dx\le\, C(t_0,\tau,R).
\eqntag
\label{eqn;L^2-bound}
\end{align*}
Subsequently, we turn into the $L^p$-estimate for $u\varphi_{4R}$.
By the above almost same argument, we eventually obtain the following differential inequality
\begin{align*}
&\frac{\d}{\d t}\int_{\R^2}(u\varphi_{4R})^p\,\d x+\left(2-3\ep\right)\int_{\R^2 }|\nabla u^{\frac p2}|^2\varphi_{4R}^p\dx
\\
\le\,&\int_{\R^2}u^{p+1}\varphi_{4R}^p\dx+2\int_{\R^2}u^{p+1}\varphi_{4R}^{p+1}\dx+C_\ep\|u\varphi_{4R}\|_p^p\left(1+\frac{1+\|\partial_t v\|_2^2}{R^2}\right)
\\&+\frac{C_\ep\|v\|_{2(p+1)}^{2(p+1)}+C\|v\|_{p+1}^{p+1}}{R^{2(p+1)}}.
\end{align*}
Now, note that $\supp\varphi_{4R}\subset\{x; |x|>4R\}$, so that
H\"older's, Gagliardo--Nirenberg's and Young's inequalities along with the $L^2$-estimate \eqref{eqn;L^2-bound} imply that
 for any $\ep>0$
\begin{align*}
\int_{\R^2}u^{p+1}\varphi_{4R}^p\dx=\,&\int_{|x|>4R}u^{p+1}\varphi_{4R}^p\dx
\\
\le\,&\left(\int_{  |x|>4R }u^2\dx\right)^{\frac12} \| (u\varphi_{4R})^{\frac p2}\|_4^2
\\
\le\,&C(t_0,\tau,R) \|(u\varphi_{4R})^{\frac p2}\|_2\|\nabla (u\varphi_{4R})^{\frac p2}\|_2
\\
\le\,&C\|u\varphi_{4R}\|_p^{\frac p2}
\left(  \||\nabla u^{\frac p2}|\varphi_{4R}^{\frac p2}\|_2+\|u^{\frac p2}|\nabla \varphi_{4R}^{\frac p2}|\|_2 \right)
\\
\le\,&\ep\int_{\R^2}|\nabla u^{\frac p2}|^2\varphi_{4R}^p\dx+C_\ep\|u\varphi_{4R}\|_p^p\left( 1+ \frac{1}{R^2} \right).
\end{align*}
By the similar way to the above estimate,
\begin{align*}
\int_{\R^2}u^{p+1}\varphi_{4R}^{p+1}\dx\le\,&\|u\varphi_{4R}\|_2\|(u\varphi_{4R} )^{\frac p2}\|_4^2
\\
\le\,&C\left(\sup_{t_0\le t \le\tau}\|u\varphi_{4R}\|_2\right) \|(u\varphi_{4R})^{\frac p2}\|_2\|\nabla (u\varphi_{4R})^{\frac p2}\|_2
\\
\le\,&C\left(\sup_{t_0\le t \le\tau}\|u\varphi_{4R}\|_2\right) \|u\varphi_R\|_p^{\frac p2}
\left(  \||\nabla u^{\frac p2}|\varphi_R^{\frac p2}\|_2+\|u^{\frac p2}|\nabla \varphi_R^{\frac p2}|\|_2 \right)
\\
\le\,&\ep\int_{\R^2}|\nabla u^{\frac p2}|^2\varphi_{4R}^p\dx+C_\ep\|u\varphi_{4R}\|_p^p\left( \|u\varphi_{4R}\|_2^2+ \frac{1}{R^2} \right),
\end{align*}
from which it follows that
\begin{align*}
\frac{\d}{\d t}\int_{\R^2}(u\varphi_{4R})^p\,\d x+\left(2-4\ep\right)\int_{\R^2 }|\nabla u^{\frac p2}|^2\varphi_{4R}^p\dx
\le\,&C_\ep\|u\varphi_{4R}\|_p^p\left(1+\frac{1+\|\partial_t v\|_2^2}{R^2}\right)+\frac{C_\ep}{R^{2(p+1)}}.
\end{align*} 
Consequently, we obtain the $L^p$-estimate for $u\varphi_{4R}$ by choosing $\ep>0$ sufficiently small and applying Gronwall's inequality that
\begin{align*}
\|u\varphi_{4R}\|_p^p\le\,&C(t_0,\tau,R,\|u(t_0)\|_p) \exp\left(\int_{t_0}^t \left(1+\frac{1+\|\partial_t v\|_2^2}{R^2}\right)\d s\right)
\le\,C(t_0,\tau,R)
\end{align*}
for $t\in[t_0,\tau]\cap(0,T)$.
The proof is complete.
\end{pr}

\begin{lem}\label{lem;bound-infty-exterior}
Suppose assumptions as in Theorem \ref{thm;global}.
For any $t_0,\tau>0$ with $t_0,\tau$, let $R_1=R_1(t_0,\tau)>0$ be taken as in Lemma~\ref{lem;L^2-exterior}.
Then, 
for all $R>R_1$, there is a constant $C(t_0,\tau,R)>0$ satisfying
\begin{align*}
\|\nabla (v\varphi_{4R})\|_\infty+\|(u\varphi_{8R})\|_\infty\le\, C(t_0,\tau,R)
\end{align*}
for $t\in[t_0,\tau]\cap(0,T)$.
\end{lem}
\begin{pr}{Lemma \ref{lem;bound-infty-exterior}}
We first show the $L^p$-estimate for $\nabla v\varphi_{4R}$ in order to show the~$L^\infty$-bound for~$\nabla(v\varphi_{4R})$.
By use of the second equation in \eqref{eqn;KS} and integrations by parts together with well-known Bochner's formula $2\nabla v \cdot \nabla\Delta v =-2|D^2 v|^2+\Delta|\nabla v|^2$,
\begin{align*}
&\frac1p\frac{\d}{\d t}\int_{  \R^2 }|\nabla v|^p\varphi_{4R}^p\dx
\\
=\,&\int_{  \R^2 }|\nabla v|^{p-2}\nabla v \cdot \nabla(\partial_t v) \varphi_{4R}^p\dx
\\
=\,&\int_{  \R^2 }|\nabla v|^{p-2}\nabla v \cdot \nabla(\Delta v -\lambda v + u) \varphi_{4R}^p  \dx
\\
=\,&-\int_{  \R^2 }|\nabla v|^{p-2}|D^2 v|^2\varphi_{4R}^p\dx +\frac12\int_{  \R^2 }|\nabla v|^{p-2} \Delta |\nabla v|^2 \varphi_{4R}^p\dx
\\
&-\lambda \int_{  \R^2 }|\nabla v|^p \varphi_{4R}^p\dx +\int_{  \R^2 } |\nabla v|^{p-2} \nabla v \cdot \nabla u \varphi_{4R}^p \dx
\\
=\,&-(p-1)\int_{  \R^2 }|\nabla v|^{p-2}|D^2 v|^2\varphi_{4R}^p\dx
-\int_{  \R^2 }|\nabla v|^{p-2} (\nabla\varphi_{4R}^p)^TD^2 v \nabla v \dx
\\
&-\lambda \int_{  \R^2 }|\nabla v|^p \varphi_{4R}^p\dx
-\int_{  \R^2 } |\nabla v|^{p-2}\Delta v u \varphi_{4R}^p\dx
\\
&-(p-2)\int_{  \R^2 }|\nabla v|^{p-4} (\nabla v)^T D^2 v \nabla v u \varphi_{4R}^p\dx-\int_{  \R^2 } |\nabla v|^{p-2}\nabla v \cdot \nabla\varphi_{4R}^p u \dx.
\end{align*}
Since $|\nabla\varphi_{4R}^p|\le\,C\varphi_{4R}^p/R$, it follows from Young's inequality that for any $\ep>0$
\begin{align*}
\left|\int_{  \R^2 }|\nabla v|^{p-2} (\nabla\varphi_{4R}^p)^TD^2 v \nabla v \dx\right|
\le\,&\frac{C}{R}\int_{  \R^2 } |\nabla v|^{p-1} |D^2 v| \varphi_{4R}^p\dx
\\
\le\,&\ep\int_{  \R^2 }|\nabla v|^{p-2} |D^2v|^2\varphi_{4R}^p\dx+\frac{C_\ep}{R^2}\int_{  \R^2 }|\nabla v|^p\varphi_{4R}^p\dx
\end{align*}
and also by use of the trace estimate $|\Delta v | \le \sqrt{2} |D^2 v|$,
\begin{align*}
&\left|\int_{  \R^2 } |\nabla v|^{p-2}\Delta v u \varphi_{4R}^p\dx\right|
\\
\le\,&\sqrt{2}\int_{  \R^2 } |\nabla v|^{p-2} |D^2 v| u \varphi_{4R}^p\dx
\\
\le\,&\ep\int_{  \R^2 }|\nabla v|^{p-2}|D^2 v|^2\varphi_{4R}^p\dx
+C_\ep\int_{  \R^2 } |\nabla v|^{p-2} u^2\varphi_{4R}^p\dx
\\
\le\,&\ep\int_{  \R^2 }|\nabla v|^{p-2}|D^2 v|^2\varphi_{4R}^p\dx+C_\ep\int_{  \R^2 } |\nabla v|^p\varphi_{4R}^p\dx
+\int_{  \R^2 }(u\varphi_{4R})^p\dx
\\
\le\,&\ep\int_{  \R^2 }|\nabla v|^{p-2}|D^2 v|^2\varphi_{4R}^p\dx+C_\ep\int_{  \R^2 } |\nabla v|^p\varphi_{4R}^p\dx
+C(t_0,\tau,R)
\end{align*}
thanks to Lemma \ref{lem;L^2-exterior}.
Similarly,
\begin{align*}
\left| (p-2)\int_{  \R^2 }|\nabla v|^{p-4} (\nabla v)^T D^2 v \nabla v u \varphi_{4R}\dx \right|
\le\,&(p-2)\int_{  \R^2 } |\nabla v|^{p-2} |D^2 v| u \varphi_{4R}^p\dx
\\
\le\,&\ep\int_{  \R^2 }|\nabla v|^{p-2}|D^2 v|^2\varphi_{4R}^p\dx
\\&+C_\ep\int_{  \R^2 } |\nabla v|^p\varphi_{4R}^p\dx
+C(t_0,\tau,R)
\end{align*}
as well as 
\begin{align*}
\left|  \int_{  \R^2 } |\nabla v|^{p-2}\nabla v \cdot \nabla\varphi_{4R}^p u \dx \right|
\le\,&\frac{C}{R}\int_{  \R^2 } |\nabla v|^{p-1} u \varphi_{4R}^p\dx
\\
\le\,& \frac{C}{R^{\frac{p}{p-1}}}\int_{  \R^2 } |\nabla v|^p \varphi_{4R}^p\dx+\int_{  \R^2 }(u\varphi_{4R})^p\dx
\\
\le\,& \frac{C}{R^{\frac{p}{p-1}}}\int_{  \R^2 } |\nabla v|^p \varphi_{4R}^p\dx+C(t_0,\tau,R)
\end{align*}
for $t\in[t_0,\tau]\cap(0,T)$.
Hence, it follows from the above computations that
\begin{align*}
&\frac1p\frac{\d}{\d t}\int_{  \R^2 }|\nabla v|^p\varphi_{4R}^p\dx+\left[(p-1)-3\ep\right]\int_{  \R^2 }|\nabla v|^{p-2}|D^2 v|^2\varphi_{4R}^p\dx
\\
\le\,&C(R)\int_{  \R^2 } |\nabla v|^p\varphi_{4R}^p\dx+C(t_0,\tau,R).
\end{align*}
Gronwall's inequality implies by choosing $6\ep=(p-1)$ that
\begin{align*}
\int_{  \R^2 }|\nabla v|^p\varphi_{4R}^p\dx\le\,&C(t_0,\tau,R),\quad
\int_{ |x|>8R }|\nabla v |^p\dx \le C(t_0,\tau,R)
\eqntag
\label{eqn;nabla-v-L^p}
\end{align*}
for $t\in[t_0,\tau]\cap(0,T)$.

Next, from the second equation of \eqref{eqn;KS} for $v$,
\begin{align*}
\partial_t(v\varphi_{4R}  )=\Delta(v\varphi_{4R} )-\lambda v\varphi_{4R} -2\nabla v\cdot\nabla\varphi_{4R} -v\Delta \varphi_{4R} +u\varphi_{4R} 
\eqntag
\label{eqn;v-eq-exterior}
\end{align*}
and by use of integral formulation 
\begin{align*}
v\varphi_{4R} =\,&e^{-\lambda t} e^{(t-t_0)\Delta}((v\varphi_{4R})(t_0) )
\\&+\int_{t_0}^te^{-(t-s)\lambda}e^{(t-s)\Delta}\left[-2\nabla v(s)\cdot\nabla\varphi_{4R} -v(s)\Delta \varphi_{4R}  +u(s)\varphi_{4R} \right]\,\d s
\end{align*}
for $t\in[t_0,\tau]\cap(0,T)$.
Since by Lemma \ref{lem;v-Lpbound}
\begin{align*}
\sup_{0<t\le\tau}\|v(t)\|_p<\infty,\quad p\in[1,\infty),
\end{align*}
it follows from Lemma \ref{lem;L^2-exterior} and \eqref{eqn;nabla-v-L^p} that
\begin{align*}
\|\nabla(v\varphi_{4R} )\|_\infty
\le\,&e^{-\lambda t}\|\nabla( v(t_0)\varphi_{4R})\|_\infty
+\frac{C}{R}\int_{t_0}^t e^{-(t-s)\lambda}(t-s)^{-\frac1{p}-\frac12}\|\varphi_{4R}\nabla v(s)\|_p\,\d s
\\
&+\frac{C}{R^2}\int_{t_0}^t e^{-(t-s)\lambda}(t-s)^{-\frac1{p}-\frac12}\| v(s)\|_p\,\d s
\\
&+C\int_{t_0}^t e^{-(t-s)\lambda}(t-s)^{-\frac1{p}-\frac12}\| (u\varphi_{4R} )(s)\|_p\,\d s
\\
\le\,&\|\nabla v(t_0)\|_\infty+\frac{C}{R}\|v(t_0)\|_\infty+\frac{C}{R}\left(\sup_{t_0\le s\le \tau}\|\varphi_{4R}\nabla v(s)\|_p\right)\Gamma\left(\frac12-\frac1{p}\right)
\\
&+\frac{C}{R^2}\left(\sup_{t_0\le s\le \tau}\| v(s)\|_p\right)\Gamma\left(\frac12-\frac1{p}\right)
\\
&+C\left(\sup_{t_0\le s\le\tau}\|(u\varphi_{4R} )(s)\|_p\right)\Gamma\left(\frac12-\frac1p\right)
\end{align*}
for $t\in[t_0,\tau]\cap(0,T)$, provided $p \in (2,\infty)$, where $\Gamma(\cdot)$ is the Gamma function on $(0,\infty)$,
which shows
\begin{align*}
\|\nabla (v\varphi_{4R} )\|_\infty\le\, C(t_0,\tau,R).
\end{align*}
Next, by the similar way to that above argument 
\begin{align*}
\partial_t(u\varphi_{8R} )=\Delta(u\varphi_{8R})+u(\Delta\varphi_{8R}+\nabla v\cdot \nabla \varphi_{8R})-2\nabla\cdot(u\nabla \varphi_{8R}),
\eqntag
\label{eqn;1st-exterior-eq}
\end{align*}
and we have by use of the integral formulation
\begin{align*}
u\varphi_{8R}=\,&e^{(t-t_0)\Delta}((u\varphi_{8R})(t_0))
\\
&+\int_{t_0}^te^{(t-s)\Delta}\left[ u(s)(\Delta\varphi_{8R}+\nabla v(s)\cdot \nabla \varphi_{8R})-2\nabla\cdot(u(s)\nabla \varphi_{8R})  \right]\,\d s,
\eqntag
\label{eqn;integral-u-exterior}
\end{align*}
so that by Lemma \ref{lem;Lp_Lq_heat}
\begin{align*}
\|(u\varphi_{8R})\|_\infty
\le\,&\|u(t_0)\|_\infty+\frac{C}{R^2}\int_{t_0}^t (t-s)^{-\frac1{2}} \|(u\varphi_{8R})(s)\|_2\,\d s
\\
&+\frac{C}{R}\int_{t_0}^t (t-s)^{-\frac1{2}} \|(u\varphi_{8R})(s)\|_p\|\nabla v(s)\|_{L^q(|x|>8R)}\,\d s
\\
&+\frac{C}{R}\int_{t_0}^t (t-s)^{-\frac1{p}-\frac12} \|(u\varphi_{8R})(s)\|_p\,\d s
\\
\le\,&\|u(t_0)\|_\infty+\frac{C}{R^2} t^{\frac12}\left(\sup_{t_0\le s\le \tau}\|(u\varphi_{8R})(s)\|_p\right)
\\
&+\frac{C}{R}t^{\frac12}\left(\sup_{t_0\le s \le\tau}\|(u\varphi_{8R})(s)\|_p\right)\left(\sup_{t_0\le s\le \tau}\|\nabla v(s)\|_{L^q(|x|>8R)}\right)
\\
&+\frac{C}{R} t^{\frac12-\frac1{p}} \left(\sup_{t_0\le s \le\tau}\|(u\varphi_{8R})(s)\|_p\right),
\end{align*}
provided that 
\begin{align*}
\frac1p+\frac1q=\frac12,\quad p,q>2.
\end{align*}
Hence, thanks to Lemma \ref{lem;L^2-exterior}, we conclude the $L^\infty$-estimate for $u\varphi_{4R}$ that
\begin{align*}
\|(u\varphi_{8R})\|_\infty\le \,C(t_0,\tau,R),
\end{align*}
which ends the proof.
\end{pr}
\vspace{5mm}
\begin{lem}\label{lem;L2smooth}
Let $(u,v)$ be the solution to \eqref{eqn;KS}.
Let $R$ be taken as in Lemma \ref{lem;L^2-exterior}.
Then, for any $t_0,\tau>0$ with $t_0<\tau$, there is a constant $C(t_0,\tau,R)>0$ such that
\begin{align*}
\|\nabla(u\varphi_{8R})\|_2&+\|\Delta(u\varphi_{8R})\|_2
+\|\nabla\Delta(u\varphi_{8R})\|_2
\\
&+\|\Delta(v \varphi_{8R})\|_2+\|\nabla\Delta(v \varphi_{8R})\|_2\le\,C(t_0,\tau,R)
\end{align*}
for $t\in[t_0,\tau]\cap(0,T)$.
\end{lem}
\begin{pr}{Lemma \ref{lem;L2smooth}}
We first show the $L^2$-estimates for $\nabla (u\varphi_{8R} )$ and use the iterating arguments to have the further regularity estimates.
By use of \eqref{eqn;1st-exterior-eq}, we have
\begin{align*}
&\frac{\d}{\d t}\int_{\R^2}|\nabla(u\varphi_{8R})|^2\dx
\\
=\,&2\int_{\R^2}\nabla(u\varphi_{8R})\cdot\nabla \partial_t \left( u \varphi_{8R}\right)\dx
\\
=\,&2\int_{\R^2}\nabla(u\varphi_{8R})\cdot \nabla \left[ \Delta(u\varphi_{8R})+u(\Delta\varphi_{8R}+\nabla v\cdot \nabla \varphi_{8R})-2\nabla\cdot(u\nabla \varphi_{8R})  \right]\dx
\\
=\,&-2\int_{\R^2}|\Delta(u\varphi_{8R})|^2\dx-2\int_{\R^2}\Delta(u\varphi_{8R}) (u\Delta\varphi_{8R}+u\nabla v\cdot \nabla \varphi_{8R})\dx
\\
&+4\int_{\R^2}\Delta(u\varphi_{8R})\nabla\cdot(u\nabla\varphi_{8R})\dx.
\end{align*}
As for the second term on the right hand side, Young's inequality together with Lemma~\ref{lem;L^2-exterior} implies that
\begin{align*}
2\left|\int_{\R^2}  \Delta(u\varphi_{8R}) u\Delta\varphi_{8R}\dx \right|
\le\,&\ep\int_{\R^2}|\Delta(u\varphi_{8R})|^2\dx+\frac{C_\ep}{R^4}\int_{\R^2}(u\varphi_{8R})^2\dx
\\
\le\,&\ep\int_{\R^2}|\Delta(u\varphi_{8R})|^2\dx+\frac{C_\ep}{R^4}\left(\sup_{t_0\le t \le\tau}\|(u\varphi_{2R})(t)\|_2\right)^2
\end{align*}
and by Lemma \ref{lem;bound-infty-exterior}
\begin{align*}
&2\left|  \int_{\R^2}  \Delta(u\varphi_{8R})u\nabla v\cdot \nabla \varphi_{8R}\dx  \right|
\\
\le\,&\frac{C}{R}\int_{\R^2}|\Delta(u\varphi_{8R})| u |\nabla v| \varphi_{8R}\dx
\\
\le\,&\ep\int_{\R^2}|\Delta(u\varphi_{8R})|^2\dx+\frac{C_\ep}{R^2} \int_{|x|>8R}u^2( |\nabla(v\varphi_{8R})|+v|\nabla\varphi_{8R}|)^2\dx
\\
\le\,&\ep\int_{\R^2}|\Delta(u\varphi_{8R})|^2\dx+\frac{C_\ep}{R^2} \left(\sup_{t_0\le t\le\tau}\|\nabla(v\varphi_{8R} ) (t)\|_\infty\right)\|u\|_{L^2(|x|>8R)}^2
\\
&+\frac{C_\ep}{R^4}\left(\sup_{t_0\le t\le\tau}\|(u\varphi_{8R})(t)\|_\infty\right)^2\|v\|_2^2
\\
\le\,&\ep\int_{\R^2}|\Delta(u\varphi_{8R})|^2\dx+\frac{C_\ep(t_0,\tau,R)}{R^2} +\frac{C_\ep(t_0,\tau,R)}{R^4} .
\end{align*}
In addition,
\begin{align*}
\left| 4\int_{\R^2}\Delta(u\varphi_{8R})\nabla\cdot(u\nabla\varphi_{8R})\dx \right|
\le\,&\ep \int_{\R^2}|\Delta(u\varphi_{8R})|^2\dx
\\
&+\frac{C_\ep}{R^2}\int_{\R^2}|\nabla(u\varphi_{8R})|^2\dx
+\frac{C\ep}{R^4}\int_{\R^2} (u\varphi_{8R})^2\dx.
\end{align*}
Therefore, collecting above computations infers
\begin{align*}
\frac{\d}{\d t}\int_{\R^2}|\nabla(u\varphi_{8R})|^2\dx+(2-3\ep)\int_{\R^2}|\Delta(u\varphi_{8R})|^2\dx
\le\,&\frac{C_\ep}{R^2}\int_{\R^2}|\nabla(u\varphi_{8R})|^2\dx+C_\ep(t_0,\tau,R).
\end{align*}
Choosing $\ep=1/3$ and applying Gronwall's inequality, we have
\begin{align*}
\|\nabla(u\varphi_{8R})\|_2^2+\int_{t_0}^t \int_{\R^2}|\Delta(u\varphi_{8R})|^2\dx\d s\le\,C(t_0,\tau,R)( \|\nabla u(t_0)\|_2^2+1) 
\eqntag
\label{eqn;nablaL2-u}
\end{align*}
for $t\in[t_0,\tau]\cap(0,T)$.
This implies from \eqref{eqn;v-eq-exterior} and Young's inequality that
\begin{align*}
&\frac{\d}{\d t}\int_{ \R^2} |\Delta (v \varphi_{8R})|^2\dx
\\
=\,&-2\int_{ \R^2}\nabla\Delta(v \varphi_{8R})\cdot\nabla
\left(\Delta(v\varphi_{8R} )-2\nabla v\cdot\nabla\varphi_{8R} -v\Delta \varphi_{8R} -\lambda v\varphi_{8R} +u\varphi_{8R} \right)\dx
\\
\le\,&-(2-\ep)\int_{ \R^2}|\nabla\Delta (v\varphi_{8R})|^2\dx+C_\ep\int_{ \R^2}\left| \nabla (\nabla v\cdot \nabla \varphi_{8R})  \right|^2\dx+C_\ep \int_{ \R^2}|\nabla (v\Delta \varphi_{8R})|^2\dx
\\
&+C_\ep\int_{ \R^2}|\nabla(v\varphi_{8R})|^2\dx+C_\ep\int_{ \R^2}|\nabla(u\varphi_{8R})|^2\dx.
\eqntag
\label{eqn;delta-v-L2}
\end{align*}
As for the second term on the right hand side of \eqref{eqn;delta-v-L2},
\begin{align*}
&C_\ep\int_{ \R^2}\left| \nabla (\nabla v\cdot \nabla \varphi_{8R})  \right|^2\dx
\\
\le\,&\frac{C_\ep}{R^2}\int_{ \R^2}|\varphi_{8R}\nabla^2v|^2\dx+\frac{C_\ep}{R^4}\int_{ \R^2} |\nabla v|^2 \varphi_{8R}^2\dx
\\
\le\,&\frac{C_\ep}{R^2}\int_{ \R^2}|\nabla^2(v\varphi_{8R})|^2\dx+\frac{C_\ep}{R^4}\int_{ \R^2}|\nabla v|^2\varphi_{8R}^2\dx+\frac{C_\ep}{R^6}\int_{ \R^2}v^2\varphi_{8R}^2\dx
\\
\le\,&\frac{C_\ep}{R^2}\int_{ \R^2}|\Delta(v\varphi_{8R})|^2\dx+\frac{C_\ep}{R^4}\int_{ \R^2}|\nabla v|^2\varphi_{2R}\dx+\frac{C_\ep}{R^6}\sup_{0\le t\le\tau}\|v(t)\|_2^2
\eqntag
\label{eqn;v-varphi}
\end{align*}
for $t\in[t_0,\tau]\cap(0,T)$,
where we use $\varphi_{8R}\nabla^2 v=\nabla^2(v\varphi_{8R})-2(\nabla v)^T \nabla \varphi_{8R}-v\nabla^2 \varphi_{8R}$.
Similarly, for the third term in \eqref{eqn;delta-v-L2} is estimated as
\begin{align*}
C_\ep \int_{ \R^2}|\nabla (v\Delta \varphi_{8R})|^2\dx\le\,&\frac{C_\ep}{R^4}\int_{ \R^2} |\nabla v|^2 \varphi_{8R}^2\dx+\frac{C_\ep}{R^6}\int_{ \R^2} v^2\dx
\\
\le\,&\frac{C_\ep}{R^4}\int_{ \R^2} |\nabla v|^2 \varphi_{2R}\dx+\frac{C_\ep}{R^6}\sup_{0\le t\le\tau}\|v(t)\|_2^2.
\end{align*}
Hence, integrating \eqref{eqn;delta-v-L2} over $(t_0,t)$ and choosing $\ep=1$  gives
\begin{align*}
\int_{ \R^2} |\Delta (v \varphi_{8R})|^2\dx\le\,C(t_0,\tau,R)
\end{align*}
for $t\in[t_0,\tau]\cap(0,T)$ after applying Gronwall's inequality
as we already know the estimate that
\begin{align*}
\int_{t_0}^t \int_{ \R^2}|\nabla v|^2\varphi_{2R}\dx\d s + \|(u\varphi_{2R})\|_2\le\,C(t_0,\tau,R)
\end{align*}
by virtue of Lemmas \ref{lem;LlogL-bound-exterior} and \ref{lem;L^2-exterior}.
Next, using \eqref{eqn;1st-exterior-eq} again, we have by Young's inequality
\begin{align*}
&\frac{\d}{\d t}\int_{ \R^2} |\Delta (u\varphi_{8R})|^2\dx
\\
=\,&2\int_{ \R^2} \Delta(u\varphi_{8R})\Delta( \partial_t(u\varphi_{8R}) )\dx
\\
=\,&-2\int_{ \R^2} \nabla\Delta(u\varphi_{8R}) \cdot \nabla \left(\left[ \Delta(u\varphi_{8R})+u(\Delta\varphi_{8R}+\nabla v\cdot \nabla \varphi_{8R})-2\nabla\cdot(u\nabla \varphi_{8R})  \right]\right)\dx
\\
\le\,&-(2-\ep)\int_{ \R^2} |\nabla\Delta(u\varphi_{8R})|^2\dx+C_\ep\int_{ \R^2}|\nabla(u\Delta\varphi_{8R})|^2\dx
+C_\ep\int_{ \R^2}|\nabla(\nabla v\cdot \nabla\varphi_{8R})|^2\dx
\\
&+C_\ep\int_{ \R^2}|\nabla\left(\nabla\cdot(u\nabla\varphi_{8R})\right)|^2\dx.
\eqntag
\label{eqn;delta-u-L2}
\end{align*}
The second term in \eqref{eqn;delta-u-L2} is estimated by Lemma \ref{lem;L^2-exterior} and \eqref{eqn;nablaL2-u} as
\begin{align*}
C_\ep\int_{ \R^2}|\nabla(u\Delta\varphi_{8R})|^2\dx\le\,&\frac{C_\ep}{R^4}\int_{ \R^2}|\varphi_{8R}\nabla u|^2\dx
+\frac{C_\ep}{R^6}\int_{ \R^2}(u\varphi_{8R})^2\dx
\\
\le\,&\frac{C_\ep}{R^4}\int_{ \R^2}|\nabla(u\varphi_{8R})|^2\dx
+\frac{C_\ep}{R^6}\int_{ \R^2}(u\varphi_{8R})^2\dx\le\, C(t_0,\tau,R)
\end{align*}
for $t\in[t_0,\tau]\cap(0,T)$.
By the same way to that of \eqref{eqn;v-varphi},
\begin{align*}
C_\ep\int_{ \R^2}|\nabla(\nabla v\cdot \nabla\varphi_{8R})|^2\dx
\le\,&\frac{C_\ep}{R^2}\int_{ \R^2}|\Delta(v\varphi_{8R})|^2\dx+\frac{C_\ep}{R^4}\int_{ \R^2}|\nabla v|^2\varphi_{2R}\dx+\frac{C_\ep}{R^6}\sup_{0\le t\le\tau}\|v(t)\|_2^2
\\
\le\,& C_\ep(t_0,\tau,R)
\end{align*}
for $t\in[t_0,\tau]\cap(0,T)$.
As for the last term in \eqref{eqn;delta-u-L2}, the analogous argument leads us to 
\begin{align*}
C_\ep\int_{ \R^2}|\nabla\left(\nabla\cdot(u\nabla\varphi_{8R})\right)|^2\dx
\le\,&\frac{C_\ep}{R^2}\int_{ \R^2}|\varphi_{8R}\nabla^2 u|^2\dx+\frac{C_\ep}{R^4}\int_{ \R^2}|\varphi_{8R}\nabla u |^2\dx
\\&+\frac{C_\ep}{R^6}\int_{ \R^2}(u\varphi_{8R})^2\dx
\\
\le\,&\frac{C_\ep}{R^2}\int_{ \R^2}|\Delta(u\varphi_{8R})|^2\dx+\frac{C_\ep}{R^4}\int_{ \R^2}|\nabla(u\varphi_{8R})|^2\dx
\\
&+\frac{C_\ep}{R^6}\int_{ \R^2}(u\varphi_{8R})^2\dx
\\
\le\,&\frac{C_\ep}{R^2}\int_{ \R^2}|\Delta(u\varphi_{8R})|^2\dx+C(t_0,\tau,R)
\end{align*}
for $t\in[t_0,\tau]\cap(0,T)$, where we utilize Lemma \ref{lem;L^2-exterior} and \eqref{eqn;nablaL2-u}.
Consequently,
Gronwall's inequality implies from \eqref{eqn;delta-u-L2} with above computations with $\ep=1$ that
\begin{align*}
\|\Delta(u\varphi_{8R})\|_2\le\, C(t_0,\tau,R).
\end{align*}
Since
\begin{align*}
&\frac{\d}{\d t}\int_{ \R^2}|\nabla\Delta(v\varphi_{8R})|^2\dx
\\
=\,&-2\int_{ \R^2}\Delta^2(v\varphi_{8R})\Delta\left[  \Delta(v\varphi_{8R} )-2\nabla v\cdot\nabla\varphi_{8R} -v\Delta \varphi_{8R} -\lambda v\varphi_{8R} +u\varphi_{8R}  \right]\dx
\\
\le\,&-\int_{ \R^2}|\nabla\Delta(v\varphi_{8R})|^2\dx+C\int_{ \R^2}|\Delta(\nabla v\cdot \nabla\varphi_{8R})|^2\dx
+C\int_{ \R^2}|\Delta(v\Delta\varphi_{8R})|^2\dx
\\
&+C\int_{ \R^2}|\Delta(v\varphi_{8R})|^2\dx+C\int_{ \R^2}|\Delta(u\varphi_{8R})|^2\dx,
\end{align*}
we eventually show from the similar way to that of above estimates that
\begin{align*}
\frac{\d}{\d t}\int_{ \R^2}|\nabla\Delta(v\varphi_{8R})|^2\dx\le\,&\frac{C}{R^2}\int_{ \R^2}|\nabla\Delta(v\varphi_{8R})|^2\dx+C(t_0,\tau,R)
\end{align*}
for $t\in[t_0,\tau]\cap(0,T)$.
As a result, the estimates for $L^2$-norm of $\nabla\Delta(v\varphi_{8R})$ is also obtained, i.e.,
\begin{align*}
\|\nabla\Delta(v\varphi_{8R})\|_2\le\,C(t_0,\tau,R).
\end{align*}
Finally, let us prove the estimate of $L^2$-norm for $\nabla\Delta(u\varphi_{8R})$. By the above similar argument,
\begin{align*}
&\frac{\d}{\d t}\int_{ \R^2} |\nabla \Delta (u\varphi_{8R})|^2\dx
\\
=\,&-2\int_{ \R^2} \Delta^2(u\varphi_{8R}) \Delta\left(\left[ \Delta(u\varphi_{8R})+u(\Delta\varphi_{8R}+\nabla v\cdot \nabla \varphi_{8R})-2\nabla\cdot(u\nabla \varphi_{8R})  \right]\right)\dx
\\
\le\,&-\int_{ \R^2} |\Delta^2(u\varphi_{8R})|^2\dx+C\int_{ \R^2}|\Delta(u\Delta\varphi_{8R})|^2\dx
+C\int_{ \R^2}|\Delta(\nabla v\cdot \nabla\varphi_{8R})|^2\dx
\\
&+C\int_{ \R^2}|\Delta\left(\nabla\cdot(u\nabla\varphi_{8R})\right)|^2\dx.
\end{align*}
This implies that
\begin{align*}
\frac{\d}{\d t}\int_{ \R^2} |\nabla \Delta (u\varphi_{8R})|^2\dx
\le\,&\frac{C}{R^2}\int_{ \R^2} |\nabla \Delta (u\varphi_{8R})|^2\dx+C(t_0,\tau,R)
\end{align*}
for $t\in[t_0,\tau]\cap(0,T)$ as we have already shown all the proper estimates to be controlled.
Therefore, we end up with the finial estimate in this claim that
\begin{align*}
\| \nabla\Delta(u\varphi_{8R})\|_2\le\,C(t_0,\tau,R)
\end{align*}
for $t\in[t_0,\tau]\cap(0,T)$,
as desired.
\end{pr}

\vspace{5mm}
\begin{lem}\label{lem;smooth-exterior}
	Let $(u,v)$ be the solution to \eqref{eqn;KS}.
	Let $R$ be taken as in Lemma \ref{lem;L^2-exterior}, and $\alpha$ be a multi-index satisfying $1\le|\alpha|\le2$.
	Then, for any $t_0,\tau>0$ with $t_0<\tau$, there is a constant $C(t_0,\tau,R,\alpha)>0$ such that
	\begin{align*}
	\| \partial_x^\alpha (u\varphi_{8R})\|_\infty\le\, C(t_0,\tau,R,\alpha)
	\end{align*}
	for $t\in[t_0,\tau]\cap(0,T)$.
\end{lem}

\begin{pr}{Lemma \ref{lem;smooth-exterior}}
We shall only show the estimate of $L^\infty$-norm for $\nabla^2(u\varphi_{8R})$ as the lower derivative estimates are obtained similarly.
Recalling \eqref{eqn;integral-u-exterior} that
\begin{align*}
u\varphi_{8R}=\,&e^{(t-t_0)\Delta}((u\varphi_{8R})(t_0))
\\
&+\int_{t_0}^te^{(t-s)\Delta}\left[ u(s)(\Delta\varphi_{8R}+\nabla v(s)\cdot \nabla \varphi_{8R})-2\nabla\cdot(u(s)\nabla \varphi_{8R})  \right]\,\d s,
\end{align*}
we have
\begin{align*}
\|\nabla^2(u\varphi_{8R})\|_\infty\le\,&\|\nabla^2(u\varphi_{8R})(t_0)\|_\infty
\\
&+\int_{t_0}^t\|\nabla^2 e^{(t-s)\Delta}\left[ u(s)(\Delta\varphi_{8R}+\nabla v(s)\cdot \nabla \varphi_{8R})-2\nabla\cdot(u(s)\nabla \varphi_{8R})  \right]\|_\infty\,\d s
\\
\le\,&\|\nabla^2(u\varphi_{8R})(t_0)\|_\infty+C\int_{t_0}^t(t-s)^{-\frac12} \|\nabla^2 (u(s)\Delta\varphi_{8R})\|_2\,\d s
\\
&+C\int_{t_0}^t(t-s)^{-\frac12} \|\nabla^2(\nabla v(s)\cdot \nabla\varphi_{8R})\|_2\,\d s
\\&+ C\int_{t_0}^t(t-s)^{-\frac12} \|\nabla^2(\nabla\cdot (u(s)\nabla\varphi_{8R})\|_2\,\d s.
\end{align*}
Noticing that $|\nabla^k \varphi_{8R}|\le\,C\varphi_{8R}/R^k$ $(k\in\mathbb{N})$ and 
\begin{align*}
\varphi_{8R}\nabla^2u=\,&\nabla^2(u\varphi_{8R})-2(\nabla u)^T \nabla \varphi_{8R}-u\nabla^2\varphi_{8R}
\end{align*}
with $\|\nabla^2 (u\varphi_{8R})\|_2=\|\Delta(u\varphi_{8R})\|_2$, we see
\begin{align*}
\|\nabla^2 (u(s)\Delta\varphi_{8R})\|_2\le\,&\frac{C}{R^2}\|\Delta(u\varphi_{8R})(s)\|_2+\frac{C}{R^3}\|\nabla(u\varphi_{8R})(s)\|_2+\frac{C}{R^4}\|(u\varphi_{8R})(s)\|_2
\\
\le\,&C(t_0,\tau,R)
\end{align*}
for $s\in[t_0,\tau]\cap(0,T)$, according to Lemma \ref{lem;L^2-exterior} and Lemma \ref{lem;L2smooth}, and hence,
\begin{align*}
C\int_{t_0}^t(t-s)^{-\frac12} \|\nabla^2 (u(s)\Delta\varphi_{8R})\|_2\,\d s\le\,C(t_0,\tau,R).
\end{align*}
Next, Lemma \ref{lem;L2smooth} with Lemma \ref{lem;v-Lpbound} also implies that
\begin{align*}
\|\nabla^2(\nabla v(s)\cdot \nabla\varphi_{8R})\|_2
\le\,&\frac{C}{R}\|\nabla\Delta(v\varphi_{8R})(s)\|_2
\\
&+\frac{C}{R^2}\|\Delta(v\varphi_{8R})(s)\|_2
+\frac{C}{R^3}\|\nabla(v\varphi_{8R})(s)\|_2+\frac{C}{R^4}\|(v\varphi_{8R})(s)\|_2
\\
\le\,&C(t_0,\tau,R)
\end{align*}
for $s\in[t_0,\tau]\cap(0,T)$, so that
\begin{align*}
C\int_{t_0}^t(t-s)^{-\frac12} \|\nabla^2(\nabla v(s)\cdot \nabla\varphi_{8R})\|_2\,\d s\le\,C(t_0,\tau,R).
\end{align*}
Subsequently, the almost same way yields from Lemma \ref{lem;L^2-exterior} and Lemma \ref{lem;L2smooth}
\begin{align*}
 \|\nabla^2(\nabla\cdot (u(s)\nabla\varphi_{8R})\|_2\le\,&
 \frac{C}{R}\|\nabla\Delta(u\varphi_{8R})(s)\|_2+\frac{C}{R^2}\|\Delta(u\varphi_{8R})(s)\|_2
 \\
 &+\frac{C}{R^3}\|\nabla(u\varphi_{8R})(s)\|_2+\frac{C}{R^4}\|(u\varphi_{8R})(s)\|_2
 \\
 \le\,&C(t_0,\tau,R)
\end{align*}
for $s\in[t_0,\tau]\cap(0,T)$, from which it follows that
\begin{align*}
C\int_{t_0}^t(t-s)^{-\frac12} \|\nabla^2(\nabla\cdot (u(s)\nabla\varphi_{8R})\|_2\,\d s\le\,C(t_0,\tau,R).
\end{align*}
Collecting above outcomes leads to the desired estimate that
\begin{align*}
\|\nabla^2(u\varphi_{8R})\|_\infty \le\, C(t_0,\tau,R).
\end{align*}
Thus, the claim is concluded.
\end{pr}

\vspace{5mm}
Thanks to the boundedness of $\partial_x^{\alpha}u$ $(1\le\alpha\le2)$ on $|x|\ge 16 R$ based on Lemma \ref{lem;smooth-exterior}, the following lemma holds true by the parabolic regularity theory.
\begin{lem}\label{lem;lower-bdd-u}
Let $(u,v)$ be the solution to \eqref{eqn;KS}.
Let $R$ be taken as in Lemma \ref{lem;L^2-exterior}.
Then, for all $t_0,\tau>0$ with $t_0<\tau$,  there exist $x_0\in\R^2$ with $|x_0|> 16R$, $\ep_0\in(0,1)$ and $\delta>0$  such that 
\begin{align*}
u(t,x) \ge \delta,\quad ~~t\in[t_0,\tau]\cap(0,T),~~|x-x_0|\le\ep_0.
\end{align*}
\end{lem}
The proof can be found in \cite[Proposition~3.12]{NaOg16}. 
Since the argument is essentially the same and the required estimates have already been established in the previous lemmas, we omit the details.
If $T<\infty$, then the constants in Lemma~\ref{lem;lower-bdd-u} can be chosen uniformly for all $\tau<T$, since the coefficients $\nabla v$ and $\Delta v$ of the (linear) first equation in \eqref{eqn;KS} for~$u$ in~$\{|x|>16R\}$ are bounded on $[t_0,T)$, whose bounds determine the constants in the parabolic regularity estimates.

\vspace{5mm}
\subsection{Regularity estimates in interior domains}\label{sect;interior}
We in this subsection give the interior a priori estimates.
For $R>0$, define the interior cut-off function $\psi_R$ by
\begin{align*}
\psi_R\equiv 1- \varphi_{8R},
\eqntag
\label{eqn;psi}
\end{align*}
where $\varphi_{8R}$ is the exterior cut-off function defined in \eqref{eqn;cut-off}. In this case, we note that $\supp \psi_R\subset \{x; |x|\le 16R\}$.
Here, let us introduce the following cut-off Lyapunov functional $\L_R(t)$ denoted by
\begin{align*}
\L_R(t):=\int_{ \R^2}u(\ln u-1) \psi_R^2\dx-\int_{ \R^2}uv\psi_R^2\dx+\frac12\int_{ \R^2}|\nabla v|^2\psi_R^2\dx+\frac\lambda2\int_{ \R^2}v^2\psi_R^2\dx.
\end{align*}
The functional $\L_R(t)$ satisfies the following identity.
\begin{lem}\label{lem;interior-lyapunov}
Let $(u,v)$ be the solution to \eqref{eqn;KS}.
Then,
\begin{align*}
\frac{\d}{\d t}\L_R(t)+&\int_{ \R^2}u|\nabla(\ln u -v)|^2\psi_R^2\dx+\int_{ \R^2}|\partial_tv|^2\psi_R^2\dx
\\
=\,&-\int_{ \R^2}u(\ln u -v)\nabla(\ln u - v)\cdot \nabla \psi_R^2\dx-\int_{ \R^2}\partial_t v \nabla v \cdot \nabla\psi_R^2\dx.
\end{align*}
\end{lem}
Lemma \ref{lem;interior-lyapunov} with $\psi_R\equiv1$ is well-established as the functional $\L(t)$ is the usual Lyapunov functional for \eqref{eqn;KS}, so that we give the brief proof here, see \eqref{eqn;Lyapunov}.

\vspace{5mm}
\begin{pr}{Lemma \ref{lem;interior-lyapunov}}
By the first equation in \eqref{eqn;KS}, it follows from integrations by parts that
\begin{align*}
\frac{\d}{\d t}\left(\int_{ \R^2}u(\ln u - 1)\psi_R^2\dx\right)
=\,&\int_{ \R^2}\partial_t u \ln u \psi_R^2\dx
\\
=\,&-\int_{ \R^2}u|\nabla \ln u|^2\psi_R^2\dx-\int_{ \R^2}u\ln u \nabla \ln u\cdot \nabla \psi_R^2\dx
\\
&+\int_{ \R^2}u\nabla \ln u \cdot \nabla v \psi_R^2\dx +\int_{ \R^2}u\ln u \nabla v \cdot \nabla \psi_R^2\dx.
\end{align*}
Next, using the first and second equations in \eqref{eqn;KS}, we have
\begin{align*}
-\frac{\d}{\d t}\int_{ \R^2}uv\psi_R^2\dx
=\,&\int_{ \R^2}u\nabla \ln u \cdot \nabla v \psi_R^2\dx+\int_{ \R^2}uv\nabla\ln u \cdot \nabla \psi_R^2\dx
\\
&-\int_{ \R^2}u|\nabla v |^2 \psi_R^2\dx-\int_{ \R^2}uv\nabla v \cdot \nabla \psi_R^2\dx
\\
&-\int_{ \R^2}|\partial_t v|^2\psi_R^2\dx-\frac12\frac{\d}{\d t}\int_{ \R^2}|\nabla v |^2\psi_R^2\dx\\
&-\int_{ \R^2}\partial_t v \nabla v\cdot \nabla \psi_R^2\dx
-\frac\lambda2\frac{\d}{\d t}\int_{ \R^2}v^2\psi_R^2\dx.
\end{align*}
Hence, combining above two identities shows that
\begin{align*}
\frac{\d}{\d t}\L_R(t)+&\int_{ \R^2}u|\nabla(\ln u -v)|^2\psi_R^2\dx+\int_{ \R^2}|\partial_tv|^2\psi_R^2\dx
\\
=\,&-\int_{ \R^2}u(\ln u -v)\nabla(\ln u - v)\cdot \nabla \psi_R^2\dx-\int_{ \R^2}\partial_t v \nabla v \cdot \nabla\psi_R^2\dx
\end{align*}
and we conclude the proof.
\end{pr}
\vspace{5mm}
\begin{lem}\label{lem;interior-bdd-lyapunov}
Suppose assumptions as in Theorem \ref{thm;global}.
Let $R$ be taken as in Lemma \ref{lem;L^2-exterior}.
Then, for any $t_0,\tau>0$ with $t_0<\tau$, there is a constant $C(t_0,\tau,R)>0$ such that
\begin{align*}
\L_R(t) + \frac{1}{2}\int_{t_0}^t \int_{ \R^2} u |\nabla(\ln u -v)|^2\psi_R^2\dx\d s+\frac12\int_{t_0}^t\int_{ \R^2}|\partial_tv|^2\psi_R^2\dx\d s\le\,&\L_R(t_0)+C(t_0,\tau,R)
\end{align*}
for $t\in[t_0,\tau]\cap(0,T)$.
\end{lem}
\begin{pr}{Lemma \ref{lem;interior-bdd-lyapunov}}
Let us recall the conclusion in Lemma \ref{lem;interior-lyapunov}, that is,
\begin{align*}
\frac{\d}{\d t}\L_R(t)+&\int_{ \R^2}u|\nabla(\ln u -v)|^2\psi_R^2\dx+\int_{ \R^2}|\partial_tv|^2\psi_R^2\dx
\\
=\,&-\int_{ \R^2}u(\ln u -v)\nabla(\ln u - v)\cdot \nabla \psi_R^2\dx-\int_{ \R^2}\partial_t v \nabla v \cdot \nabla\psi_R^2\dx.
\end{align*}
Note that $\nabla\psi_R^2 =2 \psi_R \nabla \varphi_{8R}$ from the definition of $\psi_R$ in \eqref{eqn;psi}.
On the one hand, by Young's inequality 
\begin{align*}
\left| \int_{ \R^2}u(\ln u -v)\nabla(\ln u - v)\cdot \nabla \psi_R^2\dx\right|
\le\,&\frac{C}{R}\int_{ \R^2}u|\ln u-v| |\nabla(\ln u -v)| \psi_R \varphi_{8R}\dx
\\ 
\le\,&\frac12 \int_{ \R^2}u|\nabla(\ln u -v)|^2\psi_R^2\dx
\\
&+\frac{C}{R^2}\int_{ \R^2}u (\ln u -v)^2 \varphi_{8R}^2\dx.
\end{align*}
Thanks to the exterior estimates in Subsection \ref{sect;exterior}, the last term is bounded, i.e.,
\begin{align*}
\frac{C}{R^2}\int_{ \R^2}u (\ln u -v)^2 \varphi_{8R}^2\dx\le\, C(t_0,\tau,R)
\end{align*}
for $t\in[t_0,\tau]\cap(0,T)$.
On the other hand,
\begin{align*}
\left|\int_{ \R^2}\partial_t v \nabla v \cdot \nabla\psi_R^2\dx\right|
\le\,&\frac{C}{R}\int_{ \R^2}|\partial_t v| |\nabla v| \psi_R\varphi_{8R}\dx
\\
\le\,&\frac12\int_{ \R^2}|\partial_t v|^2\psi_R^2\dx+\frac{C}{R^2}\int_{ \R^2}|\nabla v|^2 \varphi_{8R}^2\dx,
\end{align*}
whose the second term is also bounded by the exterior estimates.
Therefore, we observe
\begin{align*}
\frac{\d}{\d t}\L_R(t)+&\frac12\int_{ \R^2}u|\nabla(\ln u -v)|^2\psi_R^2\dx+\frac12\int_{ \R^2}|\partial_tv|^2\psi_R^2\dx\le\, C(t_0,\tau,R)
\end{align*}
for $t\in[t_0,\tau]\cap(0,T)$.
Finally, integrating this over $(t_0,t)$ gives the desired inequality.
\end{pr}

\vspace{5mm}
We are in a position to show the interior a priori estimate of the entropy for $u$.
\begin{prop}\label{prop;interior-bdd-entropy}
	Suppose assumptions as in Theorem \ref{thm;global}.
	Let $R$ be taken as in Lemma~\ref{lem;L^2-exterior}.
	Then, for any $t_0,\tau>0$ with $t_0<\tau$, there is a constant $C(t_0,\tau,R)>0$ such that
	\begin{align*}
	\int_{ \R^2}(u\ln u)\psi_R^2\dx \le\, C(t_0,\tau,R)
	\end{align*}
	for $t\in[t_0,\tau]\cap(0,T)$.
\end{prop}

\begin{pr}{Proposition \ref{prop;interior-bdd-entropy}}
The proof is aligned with the similar way to that of Proposition~\ref{prop;bound-energy}.
Recall the inequality \eqref{eqn;NSS} that
\begin{align*}
\int_{D} gh\ \mathrm{d}x\le \int_{D} g\log g \ \mathrm{d}x+ M \log \left(\int_{D}e^h \ \mathrm{d}x\right)-M\log M, 
\quad M:=\int_{D} g\ \mathrm{d}x,
\end{align*}
where $D\subset\R^2$ is a bounded domain. Since $\supp\psi_R\subset\{x; |x|\le16R\}$, we then have
\begin{align*}
\int_{ \R^2} uv \psi_R^2\dx=\,& \int_{|x|\le 16R} (u\psi_R) (v\psi_R)\dx
\\
\le\,&(1-\alpha)\int_{|x|\le16R} (u\psi_R)\ln(u\psi_R)\dx
\\
&+(1-\alpha)\left(\int_{|x|\le 16R}u\psi_R\dx\right)\ln\left(\int_{|x|\le16R}\exp\left[ \frac{ (v\psi_R) }{1-\alpha}\right]\dx\right)
\\
&-(1-\alpha)\left(\int_{|x|\le 16R}u\psi_R\dx\right)\ln\left(\int_{|x|\le 16R}u\psi_R\dx\right)
\end{align*}
for any $\alpha\in (0,1)$.
Here, noticing that $(v\psi_R)\in H^1_0( \{|x|\le 16R\} )$, we realize that the use of Trudinger--Moser type inequality state in Lemma \ref{lem;TM_dilechlet} is applicable, so that
\begin{align*}
\int_{|x|\le16R}\exp\left[ \frac{ (v\psi_R)}{1-\alpha}\right]\dx \le\,&C_{TM} |B_{16R}(0)| \exp\left( \frac{1}{16\pi(1-\alpha)^2}\|\nabla(v\psi_R)\|_2^2 \right)
\\
\le\,&C_{TM} |B_{16R}(0)| \exp\left( \frac{1}{16\pi(1-\alpha)^2}\||\nabla v|\psi_R\|_2^2+\frac{C}{R^2(1-\alpha)^2}\|v\|_2^2\right).
\end{align*}
Hence,
\begin{align*}
&\int_{ \R^2} uv \psi_R^2\dx
\\\le\,&(1-\alpha) \int_{|x|\le16R} (u\psi_R)\ln(u\psi_R)\dx+\frac{1}{16\pi(1-\alpha)}\left(\int_{|x|\le 16R}u\psi_R\dx\right)\left(\int_{|x|\le16R}|\nabla v|^2\psi_R^2\dx\right)
\\
&+\frac{C\|u_0\|_1}{R^2(1-\alpha)}\left(\sup_{0\le t\le\tau}\|v(t)\|_2\right)^2+(1-\alpha)\|u_0\|_1\ln(C_{TM}|B_{16R}|)+\frac1e
\\
\le\,&(1-\alpha) \int_{|x|\le16R} (u\psi_R)\ln(u\psi_R)\dx
\\
&+\frac{1}{16\pi(1-\alpha)}\left(\int_{|x|\le 16R}u\psi_R\dx\right)\left(\int_{\R^2}|\nabla v|^2\psi_R^2\dx\right)+C(R,\tau,\alpha).
\end{align*}
In particular, by the exterior estimate for $u$ 
\begin{align*}
\int_{|x|\le16R} (u\psi_R)\ln(u\psi_R)\dx=\,&\int_{|x|\le 16R}(u\ln u) \psi_R\dx
+\int_{|x|\le 16R} u(\psi_R\ln\psi_R)\dx
\\
=\,&\int_{|x|\le 16R}(u\ln u) \psi_R^2\dx+\int_{8R\le |x|\le 16R}(u\ln u) \psi_R(1-\psi_R)\dx
\\&+\int_{|x|\le 16R} u(\psi_R\ln\psi_R)\dx
\\
\le\,&\int_{\R^2}(u\ln u) \psi_R^2\dx+C(t_0,\tau,R)+ \|u_0\|_1.
\end{align*}
Moreover, Lemma \ref{lem;lower-bdd-u} together with the mass conversion implies that
\begin{align*}
\int_{|x|\le 16R}u\psi_R\dx\le\,&\int_{ \R^2}u\dx - \int_{  |x|>16R }u\dx
\\
\le\,&\int_{ \R^2}u_0\dx -\int_{ |x-x_0|\le\ep_0} u\dx
\\
\le\,&8\pi - \delta |B_{\ep_0}(x_0)|
\end{align*}
for some $\delta>0$ and $\ep_0>0$.
Consequently,
\begin{align*}
\int_{ \R^2} uv \psi_R^2\dx
\le\,&(1-\alpha)\int_{\R^2}(u\ln u) \psi_R^2\dx+\frac{8\pi-\delta|B_{\ep_0}(x_0)|}{16\pi(1-\alpha)}\int_{\R^2}|\nabla v|^2\psi_R^2\dx+C(R,\tau,\alpha)
\end{align*}
for $t\in[t_0,\tau]\cap(0,T)$, from which it follows from Lemma \ref{lem;interior-bdd-lyapunov} that
\begin{align*}
\int_{ \R^2}(u\ln u)\psi_R^2\dx=\,&\L_R(t) +\int_{ \R^2}uv\psi_R^2\dx-\frac12\int_{ \R^2}|\nabla v|^2\psi_R^2\dx
-\frac\lambda2\int_{ \R^2}v^2\psi_R^2\dx+\int_{  \R^2 }u \psi_R^2\dx
\\
\le\,&\L_R(t_0)+C(t_0,\tau,R,\alpha)+(1-\alpha)\int_{\R^2}(u\ln u) \psi_R^2\dx
\\
&-\frac12\left[ 1- \frac{8\pi-\delta|B_{\ep_0}(x_0)|}{8\pi(1-\alpha)}\right]\int_{\R^2}|\nabla v|^2\psi_R^2\dx,
\end{align*}
so that
\begin{align*}
\alpha\int_{ \R^2}(u\ln u)\psi_R^2\dx\le\,&-\frac12\left[ 1- \frac{8\pi-\delta|B_{\ep_0}(x_0)|}{8\pi(1-\alpha)}\right]\int_{\R^2}|\nabla v|^2\psi_R^2\dx+C(t_0,\tau,R,\alpha)
\end{align*}
for any $\alpha\in(0,1)$. Choosing $\alpha>0$ sufficiently small such that
\begin{align*}
1- \frac{8\pi-\delta|B_{\ep_0}(x_0)|}{8\pi(1-\alpha)}>0
\end{align*}
gives the a priori estimate for the entropy
\begin{align*}
\int_{ \R^2}(u\ln u)\psi_R^2\dx\le\, C(t_0,\tau,R)
\end{align*}
for $t\in[t_0,\tau]\cap(0,T)$. Hence, we conclude the proof.
\end{pr}

\vspace{5mm}
\subsection{Global existence of solutions}\label{sect;main}
In this subsection, we finally give the proof of Theorem~\ref{thm;global}.
To this end, we collect the a priori estimates in both exterior and interior domains.

\vspace{5mm}
\begin{pr}{Theorem \ref{thm;global}}
Let $R$ be taken as in  Lemma \ref{lem;L^2-exterior}.
By Lemma \ref{lem;LlogL-bound-exterior},
\begin{align*}
\int_{  |x|>2R } (1+u) \ln(1+u) \dx\le\, C(t_0,\tau,R).
\end{align*}
On the other hand, Lemma \ref{prop;interior-bdd-entropy} gives
\begin{align*}
\int_{ \R^2}(u\ln u)\psi_R^2\dx\le\, C(t_0,\tau,R)
\end{align*}
for $t\in[t_0,\tau]\cap(0,T)$.
Now, according to \cite[Lemma 2.3]{NaOg16},
\begin{align*}
\int_{ \Omega} (1+u)\ln(1+u)\dx\le\,2\int_{ \Omega} u |\ln u | \dx+(2\ln 2)\int_{ \Omega} u\dx,
\end{align*}
where $\Omega$ is a measurable set in $\R^2$,
 so that since $\supp\psi_R\subset \{x; |x|\le 16 R\}$
 \begin{align*}
 \int_{  |x|<8R }(1+u)\ln(1+u)\dx\le\,&2\int_{  |x|<8R }u|\ln u| \dx+(2\ln 2)\|u_0\|_1
 \\
 \le\,\,&2\int_{  \R^2 }u|\ln u| \psi_R^2\dx+(2\ln 2)\|u_0\|_1
 \\
 =\,&2\int_{  \R^2 }(u\ln u) \psi_R^2\dx-4\int_{  \R^2 }(u\ln u)_-\psi_R^2\dx+(2\ln 2)\|u_0\|_1
 \\
 \le\,&C(t_0,\tau,R) +\frac4e|B_{16R}(0)|
 \end{align*}
 for $t\in[t_0,\tau]\cap(0,T)$, where we use $x\ln x \ge -1/e$ for $x\ge0$.
 This implies that
 \begin{align*}
&\int_{  \R^2 }(1+u)\ln(1+u)\dx
\\=\,&\int_{  |x|>2R } (1+u) \ln(1+u) \dx
+\int_{  |x|\le 2R } (1+u) \ln(1+u) \dx
\\
\le\,&\int_{  |x|>2R } (1+u) \ln(1+u) \dx+\int_{  |x|< 8R } (1+u) \ln(1+u) \dx
\\
\le\,&C(t_0,\tau,R) 
 \end{align*}
  for $t\in[t_0,\tau]\cap(0,T)$.
  Hence, we have along with Proposition \ref{prop;bound-energy}
  \begin{align*}
\int_{  \R^2 }(1+u) \ln (1+u)\dx+\int_{t_0}^t \|\partial_tv\|_2^2\,\d s\le\, C(t_0,\tau,R)
  \eqntag
  \label{eqn;uniform-ulogu}
  \end{align*}
    for $t\in[t_0,\tau]\cap(0,T)$.
By the parabolic regularity argument, we obtain a uniform $L^2$-bound for~$u$
on $[t_0,\tau]\cap(0,T)$ from \eqref{eqn;uniform-ulogu} (cf. \cite{Mi13,Na-Se-Yo}), which yields further regularity estimates.
In particular, the solution~$(u,v)$ cannot blow up in finite time.
Therefore, the solution to \eqref{eqn;KS} exists globally in time.
\end{pr}

\vspace{5mm}
\noindent
{\bf Acknowledgments.}
The author is supported by JSPS Early-Career Scientists, Grant Number 25K17274 and JSPS Fellows, Grant Number 25KJ0279.
The author is also supported by the French government under the program Bourses France Excellence Japon 2025.

\vspace{5mm}
\noindent
{\bf Data Availability.} 
Our manuscript has no associated data.

\vspace{5mm}

\noindent
{\bf Conflict of interest.} 
The  author states that there is no Conflict of interest.

\end{document}